\documentclass[12pt]{amsart}
\usepackage{epsfig,amssymb,epic,eepic,color}
\usepackage[all]{xy}
\numberwithin{equation}{section}
 \usepackage[T1]{fontenc}

\theoremstyle{plain}    
\newtheorem{thm}{Theorem}[section]

\newtheorem{prop}[thm]{Proposition}
\newtheorem{lemme}[thm]{Lemma}
\newtheorem{cor}[thm]{Corollary}
\newtheorem{defn}[thm]{Definition}

\theoremstyle{definition}

\newtheorem{remarque}[thm]{Remark}
\newtheorem{remarques}[thm]{Remarks}
\newtheorem{ex}[thm]{Example}
\newtheorem{rien}[thm]{}
\newtheorem{notation}[thm]{Notation}

\numberwithin{equation}{section}
\newcommand{\be}{\begin{enumerate}}
\newcommand{\ee}{\end{enumerate}}
\newcommand{\bi}{\begin{itemize}}
\newcommand{\ei}{\end{itemize}}

\def\A{\mathbb{A}}
\def\R{\mathbb{R}}

\def\Z{\mathbb{Z}}

\def\F{\mathcal{F}}
\def\G{\mathcal{G}}

\def\cO{\mathcal{O}}
\def\S{\mathcal{S}}
\def\M{\mathcal{M}}

\def\a{\alpha}
\def\om{\omega}

\def\ga{\gamma}    
\def\Ga{\Gamma}

\def\al{\alpha}
\def\be{\beta}
 
\def\de{\delta}
\def\De{\Delta}
\def\vp{\varphi}

\def\la{\lambda}
\def\La{\Lambda}

\def\si{\sigma}
\def\Si{\Sigma}

\def\ep{\varepsilon}

\def\ds{\displaystyle}

\def\nd{\noindent}
\def\bull{\hfill$\Box$\\}
\def\proof{\nd {\bf Proof.\ }}

\def\bD{\mathbb{D}}
\def\bR{\mathbb{R}}
\def\sG{\mathscr{G}}
\def\cG{\mathcal{G}}
\def\cH{\mathcal{H}}
\def\cL{\mathcal{L}}
\def\p{\partial}

\def\cO{\mathcal{O}}
\def\bS{\mathbb{S}}
\def\cT{\mathcal{T}}
\def\cD{\mathcal{D}}

\def\x{\times}

\DeclareMathOperator{\supp}{supp}

\DeclareMathOperator{\Hom}{Hom}

\DeclareMathOperator{\Orb}{Orb}

\usepackage{hyperref}
\usepackage{mathrsfs}
\usepackage{color,soul} 
\setstcolor{red} 
\setul{0.2ex}{0.4ex} 
\usepackage[usenames,dvipsnames]{xcolor}

\newcommand{\prive}{\smallsetminus}

\newcommand{\parent}[1]{\left( #1 \right) }
\newcommand{\bparent}[1]{\bigl( #1 \bigr) }

\newcommand{\sing}[1]{\left\lbrace #1\right\rbrace }

\newcommand{\ens}[2]{\left\lbrace #1 \ \middle| \ #2\right\rbrace} 
 
\begin{document}
\today

\vskip 1cm
\begin{center}{\sc Homoclinic bifurcation in Morse-Novikov theory,

 a doubling phenomenon}
 \vskip 1cm
 
 {\sc Fran\c cois Laudenbach \& Carlos Moraga Ferr\'andiz}

\end{center}
\title{}
\author{}
\address{Universit\'e de Nantes, LMJL, UMR 6629 du CNRS, 44322 Nantes, France}
\email{francois.laudenbach@univ-nantes.fr}

\address{36, av. Camille Gu\'erin, 44000 Nantes, France}
\email{crlsmrgf@gmail.com}



\thanks{Both authors were supported by ERC Geodycon;
the second author was also supported by the Japan Society for the Promotion of Science, under the ``FY2013 Postdoctoral Fellowship (short-term) for North American and European Researchers'' program.}

\begin{abstract}  We consider a compact manifold of dimension greater than 2 and a differential form of 
degree one which is closed but non-exact. This form, viewed as a multi-valued function
has a gradient vector field with respect to any Riemannian metric. After 
S. Novikov's
work and a complement by J.-C. Sikorav, under some genericity assumptions
these data yield a complex, called today the Morse-Novikov complex. Due to the non-exactness
of the form, its gradient has a
non-trivial dynamics in contrary 
to gradients of functions. In particular, 
it is possible  that the gradient  has a homoclinic orbit. The one-form being fixed, we investigate  the 
codimension-one stratum in the space of  gradients which is formed by gradients having one simple
homoclinic orbit. Such a stratum S  breaks up into a left and a right part separated by a substratum. 
The algebraic effect on the Morse-Novikov complex of crossing S
depends on the part, left or right, which is crossed. The sudden creation of infinitely many new
heteroclinic orbits may happen. Moreover, some gradients with a simple homoclinic orbit are approached 
by gradients with a simple homoclinic orbit of double energy. These two phenomena are linked.


\end{abstract}

\maketitle
\thispagestyle{empty}
\vskip 1cm
\section{Introduction}

 \begin{rien}{\sc Morse-Novikov Theory setup.} {\rm We are given a closed connected  $n$-dimensional manifold $M$ equipped with a closed differential form of degree one
 $\al$ of type Morse, meaning that its zeroes are non-degenerate.
  In other words, the local primitives 
 $f_{loc}$ of this 1-form are Morse functions. Morse-Novikov Theory deals with the case where $\al$ is non-exact, that is
 the cohomology class $u$ of $\al$ is non-zero in $H^1(M;\R)$. The set of zeroes  of $\al$ will be denoted $Z(\al)$; 
 each zero $p\in Z(\al)$ has a Morse index $i(p)\in \{0, \ldots, n\}$. The set of  zeroes of index $k$ is noted $Z_k(\al)$.  
 
 Since the undeterminacy of $f_{loc}$ is just an additive constant, for any Riemannian metric the {\it descending} or 
 {\it negative}
  gradient $-\nabla f_{loc}$ is globally defined. 
  Such a vector field $X$ will be said an $\al${\it -gradient}.
  The zeroes of $X$ coincide with the zeroes of $\al$ and are hyperbolic. Therefore, each zero $p$ of $\al$
  has a stable manifold $W^s(p,X)$ and an unstable manifold $W^u(p,X)$. Both are manifolds which are injectively 
  immersed\footnote{When $\al$ is exact (case of Morse Theory), the stable and unstable maniolds are embedded.}
   in $M$; the unstable (resp. stable) manifold is diffeomorphic to $\R^{i(p)}$ (resp.  $\R^{n-i(p)}$). In the present article,
   we will point out some dynamical particularities of gradients of {\it multivalued Morse functions}, according to the 
   terminology introduced by S. Novikov for speaking of Morse closed 1-forms \cite{novikov}.
   

   By Kupka-Smale's theorem \cite{peixoto, palis}, 
   generically  among the $\al$-gradients
   the invariant manifolds $W^u(p, X)$ and $W^s(q,X)$ are mutually transverse for every $p,q\in Z(\al)$. Note that this 
   property is not open in general. 
   A vector field whose zeroes are hyperbolic and which fulfils this transversality property
   will be named a {\it Kupka-Smale} vector field in what follows though the classical definition is more restrictive;\footnote{ In the literature, a vector field is said to be Kupka-Smale 
   if the  periodic orbits are all hyperbolic and their center-stable and 
   center-unstable manifolds are mutually transverse.} for brevity, we shall speak of $KS$ vector fields.
  In that case, if an orbit of $X$ is a {\it connecting orbit}  going from $p$ to $q$ then,  as in Morse Theory, we have
\begin{equation}
i(p)>i(q).
\end{equation} 
In particular, as $p\neq q$ such an orbit is heteroclinic.  }
\end{rien}

\begin{rien}{\sc A key fact.} \label{deep} 
{\rm
Let us define the  $\al$-{\it length} of an $\alpha$-gradient orbit $\ell$ by 
\begin{equation}\label{eq:LengthCO}
\cL(\ell):= -\int_\ell \al.
\end{equation}
This positive number is nothing but the Riemannian \emph{energy} of $\ell$. 
It may be infinite which is the case for almost every orbit when the local primitive has no
critical points of extremal Morse index. 
Here is a {\sc key fact} which makes Morse-Novikov gradient dynamics
very special; its proof, indeed elementary, will be given in Appendix A (also \cite[Prop. 2.8]{latour}).

\begin{itemize}
\item[]{\it Assume $X$ is $KS$. Let $p\in Z_k(\al)$ and $q\in Z_{k-1}(\al)$. 
Then, for every $L>0$,
the number of connecting orbits from $p$ to $q$ whose $\al$-length is bounded by $L$ is finite.}
 \end{itemize}

 \nd This fact might be the basic observation when Novikov discovered his famous {\it Novikov ring}. 
Anyhow, it gives 
a way for reading an algebraic ``counting'' of the gradient orbits which connect two zeroes of $\al$
 when the difference of  their Morse indices is 
exactly one. We are going to sketch how this counting is made. This will be 
relevant for dynamical information.
}
\end{rien}

\begin{rien}{\sc Introduction to the Morse-Novikov complex.}{\rm

In our paper, the Morse-Novikov complex is just a tool for  \emph{encoding} a part of the dynamics of $\al$-gradients: 
namely, the 
count of orbits connecting two zeroes whose Morse indices differ by one.
We adopt a point of view which is due to J.-C. Sikorav \cite{sikorav}. 
For him, the \emph{universal} Novikov ring $\La_u$
is some \emph{completion} of the group ring $\Z[\pi_1(M,*)]$ of the fundamental group
 associated with the cohomology class $[\al]=u$
 (see Novikov Condition (\ref{eq:NovikovCond})).

The Morse-Novikov complex
is only defined when the $\al$-gradient $X$ is Kupka-Smale. As a graded module,
$NC_*(\al, X)$  in degree $k$ is the free module generated by $Z_k(\al)$ over the ring $\La_u$.
Each unstable manifold $W^u (p,X)$ is oriented arbitrarily. 
The differential $\partial^X: NC_k(\al, X)\to NC_{k-1}(\al, X)$ has the following form on a generator $p\in Z_k(\al)$:
\begin{equation}
\partial^Xp= \sum_{q\in Z_{k-1}(\al)}\left( \sum_{\ga\in \Ga_p^q}n_\ga(p,q) \ga\right) q
\end{equation}
Here, $n_\ga(p,q)$ is a relative integer, $\Ga_p^q$ denotes the set of homotopy classes of paths from $p$ to $q$.
As the unstable manifolds are oriented, and hence, the stable manifolds are co-oriented, each connecting orbit
from $p$ to $q$ carries a sign (once some conventions are fixed). The integer $n_\ga(p,q)$ is the total number of the
signs carried by the connecting orbits in the homotopy class $\ga$. If $n_\ga(p,q)\neq 0$, there is at least one connecting 
orbit from $p$ to $q$ in the class $\ga$. 
One checks that $\partial^X\circ \partial^X=0$
which justifies the term ``complex'' in the algebraic sense.

By Stokes' formula, all connecting orbits in a given
homotopy class $\ga$ have the same $\al$-length. Thus, by the {\sc key fact}, $n_\ga(p,q)$ is finite. Moreover,
 by definition of the Novikov ring $\La_u$
(see \cite{farber} or  (\ref{eq:NovikovCond})), the {\it incidence coefficient}
\begin{equation}\label{incidence}
n(p,q):=  \sum_{\ga\in \Ga_p^q}n_\ga(p,q) \ga
\end{equation}
is an element of $\La_u$, and hence, $\partial ^X$ is a morphism of graded $\La_u$-rings of degree -1. Note the following:
If $n(p,q)$ is an infinite series then there is a  sequence  $(\ga_j)$ of homotopy classes in $\Ga_p^q$
whose $\al$-lengths go to $+\infty$ and such that each $\ga_j$ contains at least one 
connecting orbit.}
\end{rien}

\begin{rien}{\sc Homoclinic bifurcation.}\label{list} 
{\rm

The complement of the set of $KS$ gradients is stratified thanks to a measure of the default to be a Kupka-Smale
vector field.
Here,  we list the strata\footnote{ At the beginning, the stratification in question is just a collection of disjoint subsets of  the space of considered gradients. Theorem \ref{thm1} gives some of them the status of genuine {\it strata}. }
 of ``codimension one'': 
\begin{itemize}
\item[1)] There is a unique
pair $(p,q)$ with $i(p)=i(q)+1$ and a unique $X$-orbit from $p$ to $q$ along which $W^u(p, X)$ 
shows a minimal transversality defect to $W^s(q,X)$.
Crossing this stratum consists of creating or cancelling a pair of orbits of opposite sign. Such a pair does not appear 
in the counting above mentioned. So, we neglect the strata of this type.
\item[2)] There is a unique pair $(p,p')$,  with  $i(p)= i(p')$ and $p\neq p'$, and a unique connecting orbit from $p$ to
$p'$. Crossing such a stratum consists of making a {\it handle slide} in the sense of Morse Theory. 
This is out of the scope of our paper.
\item[3)] There is a unique  homoclinic orbit from $p$ to itself which is {\it simple}\footnote{ 'Simple' in the multiplicity sense---see Subsection (\ref{ssec:tube}).} and {\it non-broken}\footnote{ A broken orbit is the concatenation of orbits, one being linked to the next one by some zero.}.
We are only interested in the study of this bifurcation which we call {\it homoclinic bifurcation}. 
This was also considered
  by { M. Hutchings} 
  \cite{hutchings} (see also \cite{moraga}).
\end{itemize}

If the gradient $X$ belongs to the stratum $\S$ under consideration in 3), the unique homoclinic orbit $\ell$ of $X$
forms a loop with the zero $p\in Z(\al)$ as a base point. As previously said, the $\al$-length is positive and by 
Stokes' formula, the loop $\ell$ is not homotopic to zero. Let $g$ denote its homotopy class in $\pi_1(M,p)$ and let 
denote by $\S_g\subset \S$ the stratum made of the $\al$-gradients whose homoclinic orbit belongs
 to the homotopy class
$g$. Actually, $\S$ is the disjoint union  $\{\mathop{\sqcup}\limits_g\S_g\mid g\in \pi_1(M,p)\}$. 

From now on, we are going to make some restrictions on the underlying Riemannian metric: we impose
the metric to make the $\al$-gradient {\it adapted}, meaning that
 it is linearizable at every zero $p\in Z(\al)$ with a spectrum in $\{-1,1\}$ 
(See Definition \ref{df:Adapted}). Up to some rescaling, this 
 spectrum condition means that in linearizing coordinates 
about $p$ the $\al$-gradients are \emph{radial} on both the local stable and unstable manifolds. 
Let $\F_\al$ denote the space 
of adapted $\al$-gradients. 

This constraint on the specrum (not at all generic) is more informative than a metric yielding a {\it non-resonant} spectrum.\footnote{
This brings us back to a principle that
R. Thom strongly defended in the seventies: a non-generic object is richer than one of its generic approximations
since it contains 
the information of its {\it universal unfolding}.\label{universal}} 
Our constraint, for which 
Kupka-Smale's Theorem still applies (among the adapted $\al$-gradients with given germs), 
{\it enriches} the holonomy of $\ell$
 so much  that there is a well-defined real function $\chi:\S_g\to\R$ which depends only on the linearized holonomy 
 of $\ell$ from $p$ to itself and is continuous. This function will be constructed in Section \ref{sect:Sg})
 and will plays the main r\^ole in our paper.
We name  $\chi$ the {\it character function}. We have the following statement.

  \begin{thm}\label{thm1}${}$
  
   \begin{enumerate}
  \item 
  The stratum $\S_g$ is a 
  codimension-one, co-oriented submanifold of $\F_\al$ of class $C^\infty$.
  \item Assume $n>2$ and $\S_g\neq\emptyset$. 
  Then,  the vanishing locus $\S^0_g$ of $\chi$ is a non-empty
  co-oriented codimension-one submanifold of class $C^\infty$ in $\S_g$ meeting 
  each of its connected components.
 \end{enumerate}
  \end{thm}

\begin{figure}[ht]
\includegraphics[width=9cm, height=4cm]{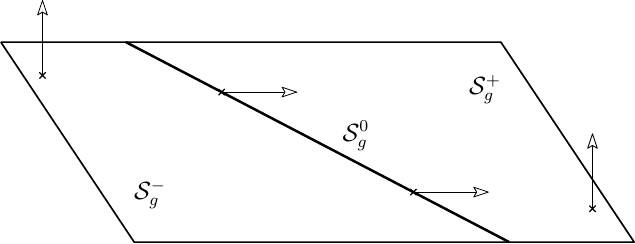}
\caption{Local situation in $\F_\al$ of every connected component of $\S_g$ near $\S_g^0$.  
The arrows indicate the co-orientation of $\S_g$ in $\F_\al$ and $\S_g^0$ in $\S_g$.}
\label{fig:Sg0dsSg}
\end{figure}
Set $\S_g^+: =\{X\in \S_g\mid \chi(X)>0\}$ and $\S_g^-: =\{X\in \S_g\mid \chi(X)<0\}$.
We shall see that crossing  $\S_g$ positively through $\S_g^+ $ or $\S_g^- $ changes the Morse-Novikov complex
in a completely different manner.
}
\end{rien}
\begin{rien}{\sc Bifurcation by crossing $\S_g$.} {\rm

Since $ \S_g$ is co-oriented, we can study the generic 
 one-parameter families $ 
 (X_s)_{s\in {\cO p(0)}}$ which intersect 
$\S_g$ \emph{positively} at $X_0$; here, we use Gromov's notation: ${\cO p(0)}$
 stands for an open interval which contains 0 and whose size is chosen as small as desired.\footnote{More generally, if
 $A$ is a closed subset of $B$, $\mathcal Op(A)$ stands for an open
 neighborhood of $A$ in $B$ which is not specified.} 
 In other words, the path $(X_s)_{s\in {\cO p(0)}}$ has to be thought as a germ of path at $X_0$.

We recall $p$, the zero of $\al$ which is involved in the homoclinic orbit of $X_0$.
Let  $q\in Z(\al)$ be any zero whose Morse index satisfies $i(q)=i(p)-1$.

The next theorem requires some genericity assumption, namely the property for $X\in \S_g$ to be
 {\it almost Kupka-Smale}
(Definition \ref{almost}). 
The corresponding residual set is noted $\S_{g,\infty}$.
Let $X_0\in \S_{g,\infty}$. Take any $L>0$. By Proposition \ref{g-residual}, for every $s\neq 0$ close enough to 0
the vector field $X_s$ is {\it Kupka-Smale up to $L$}, that is, the algebraic number of
 connecting orbits from $p$ to $q$ with
$\al$-length smaller than $L$ is finite  and locally constant. Denote it by  $n(p,q)_L^-$ or  $n(p,q)_L^+$ 
depending on the sign of $s$; these numbers are called {\it truncated incidence coefficients}.  

 Theorem \ref{thm:selfslideSimplif} states how these truncated incidence coefficients change through crossing $\S_g$.\footnote{
This makes more  precise the statement of 
\cite[Prop. 2.2.36]{moragaTesis} whose proof was inaccurate.}
More explanation of the formulas will be given 
just after the statement.

\begin{thm}\label{thm:selfslideSimplif}${}$

Let $(X_s)_{s\in {\cO p(0)}}$ be a path crossing $\S_g$ positively at time $s=0$. 
If $X_0\in \S_{g,\infty}$ 
 the following holds
for every $L>0$: 
\begin{enumerate}
\item when $X_0$ belongs to $\S_g^-$, then  $n(p,q)_L^+ =(1+g)\cdot n(p,q)_L^-$\quad{\rm (mod. $L$)},
\item when $X_0$ belongs to $\S_g^+$, then $n(p,q)_L^+ =(1+g+ g^2+ \cdots)\cdot n(p,q)_L^-$\quad{\rm (mod. $L$)}.
\end{enumerate}
\end{thm}

Of course, in order to keep the squared differential equal to zero there are similar formulas 
for the change of the incidence $n(q',p)_L^\pm$ when the Morse indices satisfy $i(q')=i(p)+1$. More precisely, we have:
{\it \begin{itemize}
\item[(3)] when $X_0\in \S_g^-$, then $n(q',p)_L^+= n(q',p)_L^-\cdot (1-g+g^2- g^3+ \cdots)$\quad{\rm (mod. $L$)},
\item[(4)] when $X_0\in \S_g^+$, then $n(q',p)_L^+= n(q',p)_L^-\cdot (1-g)$\quad{\rm (mod. $L$)}.
\end{itemize}}
\smallskip

The explanation for the product in Formulas (1) and (2) goes as follows (and similarly for (3) and (4)).  
On the one hand, recall
Formula (\ref{incidence}); namely,
 $n(p,q)= \sum_\ga n_\ga(p,q)\ga$
where $\ga\in \Ga_p^q$ is a homotopy class of paths from $p$ to $q$ and $n_\ga(p,q)$ is a relative integer
which gives the total algebraic number of connecting orbits in that class.
On the other hand $g$ is a homotopy class of loops based at $p$. 
So, the concatenation $g^j\cdot \ga$ makes sense and yields a new element 
in $\Ga_p^q$. This product is distributive
with respect to the sum.\\

}
\end{rien}
\begin{remarque} 

We have examples described in \cite{l-m} where the homoclinic bifurcation is isolated. In that case, 
truncation by $L$ becomes useless. 
There are a pair $p, q$ of zeroes with $i(p)=i(q)+1$ and a one-parameter family of adapted $\al$-gradients 
$\left(X_s\right)_{s\in [-\ep,+\ep]}$ crossing $\S_g$ positively such that:
\begin{itemize}
\item[-] for every $s<0$ there is a unique heteroclinic orbit $\ell $ from $p$ to $q$; 
\item[-] for every $s>0$ there are infinitely many heteroclinic orbits from $p$ to $q$.
\end{itemize}
This infinity which appears at once
 reads as follows: one heteroclinic orbit in each homotopy class $g^j\cdot [\ell]$, for $j= 0, 1, 2, ...$; here, [-]
stands for the homotopy  class.\\
\end{remarque}

\begin{rien}{\sc The doubling phenomenon.} {\rm

This relates the strata $\S_g$ and $\S_{g^2}$. Look for instance at the above-mentioned example and consider 
a small generic loop in $\F_\al$ going around the codimension-two stratum $\S_g^0$, beginning by crossing 
$\S_g^-$ positively and returning by crossing $\S_g^+$ negatively. Then the cumulative factor after one turn
 is $(1-g^2)$ 
 while it should be equal to 1; this is a contradiction. The next theorem solves this contradiction.
 We emphasize that its statement does not need the presence of any other zero of $\al$ than the base point
 of the homoclinic orbit.

\begin{thm}\label{thm:Doubling}  Again assume $n>2$ and $\S_g\neq \emptyset$.
 Then, there exists a codimension-two stratum $\S_g^{0,0}$ in $\S_g$, contained in $\S_g^0$, such that 
$\S_g^0\smallsetminus \S_g^{0,0}$ adheres to $\S_{g^2}$ as a boundary\footnote{This does not contradict the fact that by their very definition
$\S_g$ and $\S_{g^2}$ are disjoint.} of class $C^1$, 
 more precisely of its positive part $\S_{g^2}^+$. 
\end{thm}

In dynamical language, if $X$ is an adapted $\al$-gradient in $\S_g^0$ (that is, $\chi(X)= 0$)
whose homoclinic orbit is $\ell$, then $X$
can be approximated by an adapted $\al$-gradient $X'$  having a unique homoclinic orbit 
$\ell'$ turning twice around $\ell$. In particular,  $[\ell']= [\ell]^2$ in $\pi_1(M,p)$.

The precise definition of $\S_g^{0,0}$  yields a decomposition 
$\S_g^0\prive\S_g^{0,0}=\S_g^{0,-}\,\sqcup\,\S_g^{0,+}$ (see Definition \ref{decompositionS^0}). 
Locally along $\S_g^0\smallsetminus\S_g^{0,0}$, the stratum $\S_{g^2}$ approaches from  one side of 
 $\S_g$ only. As a matter of fact, $\S_{g^2}$ approaches $\S_g^{0,+}$  (resp. $\S_g^{0,-}$) 
  from the positive (resp. negative side) of the co-oriented stratum $\S_g$.
 A precise statement about the latter facts is given in Theorem \ref{thm:DoublingRefined} and 
 may be illustrated by Figure \ref{figure2}.
\begin{figure}[h] 
\includegraphics[width=9cm,height=4cm]{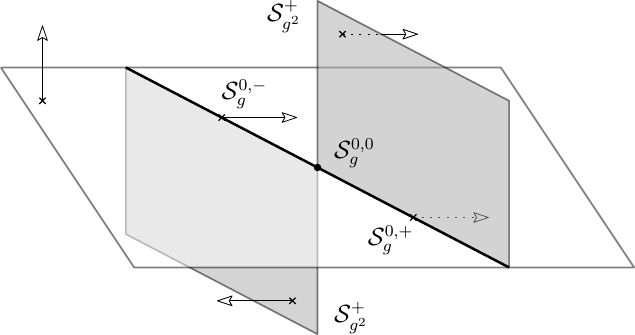}
\caption{The stratum $\S_g^0\prive\S_g^{0,0}$ as a boundary of $\S_{g^2}$ which is grayed.}\label{fig:Sg0Detail}\label{figure2}

\end{figure}
}
\end{rien}

\begin{remarques}${}$

\nd 1) Our doubling phenomenon evokes 
the period doubling bifurcation, also 
called Andronov-Hopf's bifurcation. In this aim, it would be good to know that crossing $\S_g$
creates (or destroys) a periodic orbit in the free homotopy class $g$. Shilnikov's theorem \cite{shilnikov}
deals with this question. Unfortunately, it is not applicable here because we are going to use
very non-generic Morse charts (only $+1$ and $-1$ as eigenvalues of the Hessian at critical points).
In counterpart,
 such charts offer very nice advantages.\\
 
 \nd 2) We were asked the question whether our results depend on the assumption that the
 $\al$-gradient $X$ is {\it adapted} 
 in the sense of Definition \ref{df:Adapted}. Most probably,
 if this assumption is not fulfilled, there is no way to define something which gives so much information as
 the {\it character function}. In this case, all of the three above-stated
 theorems  disappear.\footnote{ This confirms Thom's principle which we have mentioned in Footnote \ref{universal}.}\\
 
 \nd 3) We were also suggested to find a more general setting for our results. For instance, one could fix 
 a Riemannian metric and look
 at vector fields with given hyperbolic zeroes. In this setting, the key fact \ref{deep} still holds true if the 
 $\al$-length is replaced with the Riemannian energy. Unfortunately, there is no natural stratification
 of the complement of Kupka-Smale vector fields. Namely, there is no ``grading'' of the homoclinic orbits
 based at a zero of the considered vector fied. For instance, one  faces, in general, a defect of equicontinuity of sequences 
 of the homoclinic orbits in a given homotopy class. Therefore,  we do not see any natural generalization of our study.\\
 
 \nd {\bf Acknowledgement.} The first author is grateful to Yulij Ilyashenko for fruitful discussions
 on the occasion of the Conference ``Topological methods in dynamics and related topics'', Nizhny Novgorod, 2019;
 and to Slava Grines who invited him. We also thank Dirk Schuetz and Claude Viterbo for pertinent comments.
  \end{remarques}

\section{\sc Homoclinic bifurcation, orientation and character}\label{sect:Sg}

We are focusing on homoclinic bifurcations even though some of the statements
 hold true for other bifurcations (see list in Subsection \ref{list}). We consider an $\al$-gradient $X$ 
 with a simple homoclinic 
 orbit $\ell$ based at some zero $p\in Z(\al)$ in the homotopy class $g\in \pi_1(M,p)$.
 In this section, we are going to show that if the Morse coordinates about $p$ are {\it simple} in the sense 
 of Definition \ref{df:Adapted} then the Morse model $\M_p$ in these coordinates
  allows us to enrich the holonomy of $\ell$ with some specific information. 
 From this latter we deduce the {\it character function} which is the key new tool of the paper. Finally, we prove 
 Theorem \ref{thm1}.

\begin{defn}\label{df:Adapted}${}$

\nd {\rm 1)} For each zero $p$ of $\al$ of index $i(p)=i$, 
 \emph{simple Morse coordinates}
about $p$ are coordinates where the form 
$\al$ is equal to the differential of the standard quadratic form 
$$Q_i:= \frac 12\left[-x_1^2+\ldots-x_i^2+x_{i+1}^2+\ldots +x_n^2\right].
$$

\nd {\rm 2)} An $\al$-gradient $X$  is said to be \emph{adapted}  if 
 for every $p\in Z(\al)$  there are simple Morse coordinates about $p$ such that 
 $X$  coincides with the \emph{standard} descending gradient $X_{i}$
 of $Q_i$, where $i=i(p)$, that is:
  $$X_{i}=\sum_1^i x_k\partial_{x_k}-\sum_{i+1}^n  x_k\partial_{x_k}.$$
  Such  Morse coordinates are also said to be \emph{adapted} to $X$.
 \end{defn}
The property for  $X$ to be adapted 
depends only on the germ of $X$ near $Z(\al)$.
 We recall that for simplicity we fix the germ of adapted $\al$-gradients once and for all at every zero of $\al$;
 the set of such adapted $\al$-gradients  noted $\F_\al$.

  \begin{remarque}\label{alexander}
  The simplicial group $\G:=\text{Diff}(Q_i)$ of germs of diffeomorphisms of
  $(\R^n,0)$ preserving $Q_i$ retracts by deformation to $O(i,n-i)$, the linear group of isometries of $Q_i$. 
  Indeed, if 
  $\vp\in \G$,  the Alexander isotopy $\vp_t: x\mapsto \frac 1t\vp(tx)$ is made of elements in $\G$
  for every $t\in (0,1]$ and tends to the derivative $\vp'(0)x$ as $t $ goes to 0.
  Moreover, $O(i,n-i)$ retracts by deformation to its maximal compact subgroup $O(i)\times O(n-i)$
  which is the isometry group of the pair $(Q_i,X_i)$. 
  
  As a consequence, the space of germs 
  of adapted $\al$-gradients is made of a unique element up to the action of ${\rm Diff}(M\,\text{rel}\, Z(\al))$. For this reason, with 
  no loss of generallity, we may fix the germ at $p$ of all the considered adapted $\al$-gradients in what follows
  for every $p\in Z(\al)$. This choice will be done for all bifurcation families in Sections \ref{section3}  and \ref{section4}.
 \end{remarque}
 
 \subsection{\sc Morse model.} \label{ssec:morse}
 
Given $p\in Z(\al)$ of Morse index $i$,
 a {\it Morse model} $\M_p\subset M$ with positive  parameters $(\de,\de^*)$ (which we do not 
 make explicit in the notation) 
 is diffeomorphic to the subset of  
 $\R^i\times\R^{n-i}$ made of pairs $(x^-,x^+)$ such that $Q_i(x^-,x^+)\in [-\de^*,+\de^*]$, 
 $\vert x^-\vert^2\vert x^+\vert^2\leq \de\de^*$ and $\al\vert_{\M_p}= dQ_i$.  The bottom of $\M_p$, that is its intersection with  $\{Q_i=-\de^*\}$
 is denoted by $\partial^-\M_p$; { similarly, 
 the top is denoted by} $\partial^+\M_p$.
 The rest of the boundary of $\M_p$ is 
 denoted by $\partial^\ell \M_p$ and $X_i$ is tangent to it. Note that:
 \begin{itemize}
  \item the group $G$ preserves $\M_p$ for every parameters $(\de,\de^*)$;
  \item the set of Morse models, as compact subsets of $\R^n$, is contractible.
  \end{itemize}
 The flow of $X_i$ is denoted by $(X_i^t)_{t\in \R}$.
 The {\it local unstable} (resp. {\it local stable}) manifold is formed by the points
 $x\in \M_p$ whose negative (resp. positive) flow line
 $X_i^t(x)$ goes to $p$ when $t$ goes to $-\infty$ (resp. $+\infty$) without getting out of $\M_p$.
 Denote by $\Si^-$ the $(i-1)$-sphere which is formed by the points  in the bottom of $\M_p$
 which belong to   $W_{loc}^u(p,X_i)$; that is
 \begin{equation}
 \Si^-=\{(x^-,x^+)\mid  |x^-|^2=2\de^*, |x^+|= 0\}.
 \end{equation}
 This is called the {\it attaching sphere}.
  Similarly, $\Si^+$ denotes the {\it co-sphere}, the $(n-i-1)$-sphere which is contained
  in the top of $\M_p$ and made of points belonging to $W_{loc}^s(p,X_i)$, that is
  \begin{equation}
 \Si^+=\{(x^-,x^+)\mid |x^-|= 0, |x^+|^2=2\de^*\}.
 \end{equation}
  We will use the two projections associated with these coordinates:
 \begin{equation}
 \pi^+:\partial^+\M_p\to \Si^+ \quad\text{and} \quad  \pi^-:\partial^-\M_p\to \Si^-.
 \end{equation}

\subsection{\sc Simple homoclinic orbit, tube and orientation.} \label{ssec:tube} 

Let $X$ be an  adapted $\al$-gradient and let $\M_p$ be a Morse model adapted to $X$ about $p\in Z(\al)$.
 A homoclinic orbit $\ell$  of $X$ based at $p$
 is said to be {\it simple} 
 when at any point $m\in \ell$ the span $T_mW^u(p,X)+T_mW^s(p,X)$ is of codimension one in $T_mM$. When 
 $X$ is said to have a {\it unique} homoclinic orbit it will be meant that this orbit is simple, that is unique with multiplicity.
 
In this setting, denote by $\underline\ell$ the closure of $\ell\smallsetminus \M_p$; it will be named the
 {\it restricted} homoclinic orbit. The end points of $\underline \ell$
 are denoted respectively $a^-\in \Sigma^-$ 
and $a^+\in \Si^+$. We also introduce a {\it compact tube}
 $T$ around $\underline \ell$ made of   $X$-trajectories 
from $\partial^-\M_p$ to $\partial^+\M_p$. 
As $\ell$ is simple,  if the tube is small enough there are coordinates on $T$ that we note
$(x, y, v, z)\in \R^{i-1}\times  \R^{n-i-1}\times[-1,1]\times[0,1]$ with the following properties:
\begin{itemize}
\item $X$ { is positively colinear to } $\partial_z$, 
\item $\{z=0\}= T\cap \partial^-\M_p$ and  $\{z=1\}= T\cap \partial^+\M_p$;
\item $T\cap\Si^-=\{y=0, v=0, z=0\}$ and $T\cap \Si^+=\{x=0, v=0, z=1\}$;
\item $\underline\ell=\{x=0,y=0,v=0\}$;
\item the frame $(\partial_x, \partial_y, \partial_v)$ is tangent to the leaves of $\al$.
\end{itemize}
In what follows, $\{z=0\}$ (resp. $\{z=1\}$) will stand for $T\cap \partial^-\M_p$ (resp. $T\cap \partial^+\M_p$).\\

Orient the unstable $W^u(p, X)$. 
Thus, the stable manifold $W^s(p, X)$ is co-oriented. Therefore, we can choose the coordinate $v$ in the tube 
{ so that, 
for every $z_0\in [0,1]$, the following holds:}
\begin{equation}\label{orientation1}
\partial_v\wedge \text{or}\left(W^u(p, X)\cap\{z=z_0\}\right)= \text{co-or}(W^s(p,X)).
\end{equation}
If the orientation of $W^u(p, X)$  is changed, then the co-orientation of $W^s(p,X)$ 
is  also changed and the above equation shows that the positive direction 
of $v$ remains unchanged. 

\begin{remarque}\label{rem:holonomie} It is important to notice that (\ref{orientation1})
tells us nothing about the holonomy  along $\underline\ell$ of the foliation defined by $X$
(see the next subsection).
Therefore, for a given $\partial_v\in T_{a^+}(\partial^+\M_p)$, the tangent vector 
$\partial_v\in T_{a^-}(\partial^-\M_p)$ may have any position not contained in the hyperplane $\bR\sing{\p_x,\p_y}$,
depending on $X$.\\

\end{remarque}

\subsection{\sc Holonomy and perturbed holonomy.}\label{perturbed}
The foliation of $M\smallsetminus Z(\al)$ by the orbits of $X$ together its two transversals 
$\partial ^\pm\M_p$
defines a {\it holonomy} diffeomorphism
$H_X: \mathcal N^-_X\to \mathcal N^+_X$, \hfill\break that is:
\begin{itemize}
\item  $\mathcal N^-_X$  is an open connected neighborhood of $\{z=0\}$ in $\partial ^-\M_p$; 
\item $\{z=1\}\subset \mathcal N^+_X\subset \partial^+\M_p$;
\item the restriction of $H_X$ to $\{z=0\}$ is defined by $(x,y,v,0)\mapsto (x,y,v,1)$;
\item for every $a\in \mathcal N^-_X$, the image $H_X(a)$ belongs to the $X$-orbit of $a$.
\end{itemize}
By the connectedness of $\mathcal N^-_X$, the time of the flow $X^t$ for going from 
$a$ to $H_X(a)$ is continuous and hence smooth.

The existence of such holonomy diffeomorphism is an open property with respect to $X$.
More precisely, if $X'$ is a close enough approximation of $X$ in the $C^1$-topology, 
there is a {\it perturbed holonomy} diffeomorphism 
$H_{X'}$ from an open neighborhood 
 $ \mathcal N^-_{ X'}$ 
 of $\{z=0\}$ in $\partial^-\M_p$
to  an open neighborhood  $\mathcal N^+_{X'}$ of $\{z=1\}$ in $\partial^+\M_p$.

\begin{remarque} \label{batons}
It makes sense to speak of $H_{X'}^{-1}(\Si^+)\cap \{z=0\} $. It is an $(n-i-1)$-disc 
$C^1$-close to the $y$-axis in $\{z=0\}$. 
Similarly, it makes sense to speak of $H_{X'}(\Si^-)\cap \{z=1\} $. It is an $(i-1)$-disc 
close to the $x$-axis in $\{z=1\}$. 
\end{remarque}

We now state and prove the first item of Theorem \ref {thm1}. Let $p\in Z(\al)$. 
For every $g\in \pi_1(M,p)$, we consider $\S_g\subset \F_\al$, the set of adapted
$\al$-gradients which have a unique homoclinic orbit
 forming a loop based at $p$ in the homotopy class $g$; the existence of a broken homoclinic orbit is excluded 
 from $\S_g$. Recall $u$, the cohomology class of 
 the closed form $\al$; if the evaluation $u(g)$ is non-negative then $\S_g$ is empty.

\begin{prop} \label{item1} For every $g\in \pi_1(M,p)$,  
 the subset  $\S_g\subset \F_\al$ is a $C^\infty$ codimension-one submanifold of $\F_\al$, that is 
 $\S_g$ is locally 
 defined by a  regular real valued equation.
 This stratum has a canonical co-orientation.
\end{prop}

\proof 
Let $ X_0$ be any point in $\S_g$. Let $\ell$ denote
   the homoclinic orbit which forms a loop whose class $g\in \pi_1(M,p)$. 
    We intend to find a regular real valued equation for $\S_g$ near $ X_0$.
  From Remark \ref{alexander} we have the two following properties:
   \begin{itemize}
     \item the local $C^\infty$ stability  near  $p$
of the adapted $\al$-gradients ;
    \item the acyclicity of the space of Morse models adapted to $X_0$ near $p$.   
   \end{itemize}
  Therefore, the action of  the group $\text{Diff} (M)$ on $\F_\al$ reduces us to consider a local {\it slice} 
  $S\subset \F_\al$
  for this action and to look for the smoothness of $\S_g\cap S$. Namely, choose a Morse model $\M_p$
  adapted to $X_0$ and define $S:=\{X\in \F_\al \mid X=X_0 \ \text{in}\  \M_p\}$. Thanks to the stability
  property above-mentioned,
   this $S$ is indeed a local slice for the action of $\text{Diff} (M)$.
   
   We use the tube $T$ and its coordinates  as introduced in Subsection \ref{ssec:tube}. 
   The Implicit Function Theorem allows us to follow continuously, for $X$ close to $X_0$,
   a connected component $D(X)$ of $W^u(p,X)\cap \{z=1\}$  which coincides with 
   $\{y=0,v=0, z=1\}$ when $X= X_0$. Let 
   \begin{equation}
   p_v:\sing{z=1}\to\sing{ v=0, z=1}
   \end{equation}
denote the projection parallel to $\partial_v$ onto the $(x,y)$-space.
The image  $p_v(D(X))$
 is transverse to $\Si^+$. 
 The intersection is a point $a^+(X)$ which depends $C^\infty $ on $X$. Let $b(X)$
   be the point of $D(X)$ which has the same coordinates as $a^+(X)$ except the last coordinate $v$.
    Thus, the desired equation is 
   \begin{equation}\label{eq-homoclinic}
   v(b(X))=0\,.
   \end{equation}
   
   This is clearly a $C^\infty $ equation. For proving this equation is \emph{regular} it is sufficient to exhibit 
   a one-parameter family $(X_s)_{s\in\cO p(0)}$ passing through $X_0$ and satisfying the following inequality: 
   \begin{equation}\label{trans}
   \partial_s\left(v(b(X_s)\right)_{\vert s=0}> 0.
   \end{equation}
   This is easy  to perform by taking 
   $$
   X_s=X_0 + s\,g(x,y, z,v)\,\partial_v
   $$ 
    where  $g$ is a small non-negative, supported in the interior of the tube $T$ and has a positive integral along 
 $\underline\ell$\,.
  Let us check that such  $(X_s)$ fulfils (\ref{trans}). Indeed, at every point of the tube $T$ the vector $X_s$ is in the
  span$\{\partial_v, \partial_z\}$. Therefore, if $H_s$ denotes the perturbed holonomy along $\underline \ell$
  of the flow of $X_s$ then
  $b(X_s)=H_s(a^-)$ and $v(b(X_s))=s\int_0^1g\, dz$. Hence, (\ref{trans}) is fulfilled.
   
   Right after (\ref{orientation1}) we noticed  that the positive direction of $v$ does not depend on
   the chosen orientation of the unstable manifolds. Therefore, 
    (\ref{trans}) defines a canonical co-orientation of $\S_g$.
    
    What we have done is not sufficient for proving the statement. Equation (\ref{eq-homoclinic}) only solves the 
    question of existence of a homoclinic orbit at $p$ close to $\ell$. We still have to prove that $\S_g$
    does not accumulate to itself\footnote{as does a leaf of the irrational linear foliation on the 2-torus.} near $X_0$.
      More precisely, it does not exist a sequence $X_k\in \S_g$ converging to $X_0$ such that  
      \begin{equation}\label{hyp-v}
      v(b(X_k))\neq 0\ \text{ for every }k.
      \end{equation}
    Assume such a sequence exists. Let $\ell_k$ be the  unique  homoclinic orbit of $X_k$ based at $p$.
    As the sequence $\left(X_k\right)$ is close to $X_0$, the $C^0$-norm of $X_k$ is uniformly bounded
         and then the family $\left(\ell_k\right)$ is equicontinuous.  
         By Ascoli's Theorem, there is  a sub-sequence converging $C^0$ to some line $\ell_\infty$ from $p$
         to $p$. 
         
         Let us show that $\ell_\infty$ is a (possibly broken) homoclinic orbit of $X_0$.
    Indeed, let $(x_k)$ be a sequence of points with $x_k\in \ell_k$ converging to $x_\infty\in \ell_\infty$. Let $X_k^t$
    be the flow of $X_k$. For every given $t$, the sequence $\left(X_k^t(x_k)\right)$ converges to $X_0^t(x_\infty)$
    since $X_k$ tends to $X_0$ as $k$ goes to $+\infty$. Therefore, the piece of $\ell_\infty$ between $x_\infty$
    and $X_0^t(x_\infty)$ is contained in the $X_0$-orbit of $x_\infty.$
      
      If $x_k$ is close to $\ell$, then the restricted homoclinic orbit $\underline{\ell}_k$ lies in the    
      tube $T$, and hence $v(b(X_k))=0$, contradicting our assumption (\ref{hyp-v}). Therefore, there exists
    a small tube $T'$ around $\ell$  such that $\ell_k$ avoids $T'$ for every $k$ and the $C^0$-limit 
         $\ell_\infty$ as well. Finally, we get two distinct homoclinic orbits of $X_0$ based at $p$,
         one of them  being possibly broken. This is excluded by the very definition of $\S_g$.  \bull
         
   \begin{defn}\label{df:PositTrans}
  Let  $(X_s)_{s\in \cO p(0)}$ be  a one-parameter family 
 of  adapted $\al$-gradients with $X_0\in \S_g$\,.  
 This family   is said to be 
\emph{positively transverse}  to the stratum $S_g$  if it satisfies \rm{(\ref{trans})}.
 \end{defn} 
 
Let $(X_s)_{s\in \cO p(0)}$ be such a one-parameter family
 and let $H_s$ be the perturbed holonomy along $\underline \ell$ of the flow of $X_s$.
Below, we use the coordinates $(x,y,v)$ both in $\{z=0\}$ and in $\{z=1\}$.

For further use, 
 we are interested in the local solution $x_s$ of the equation
\begin{equation}\label{x_s}
(x\circ H_s)(x,0,0)=0
\end{equation} 
which is equal to $0$ when $s=0$. In this equation, the unknown is the unique point of the $x$-space in $\{z=0\}$
whose image through $H_s$ is in the $y$-space of $\{z=1\}$.
 And similarly, we consider the solution $y_s$ of the equation
\begin{equation}\label{s_s}
(y\circ H_s^{-1})(0,y,0)=0
\end{equation} 
which is equal to $0$ when $s=0$.  

\begin{lemme}\label{velocity} With the above data and notations, 
the following equality holds:
\begin{equation}
\p_s \left( v\circ H_s\right)(x_s,0,0)_{\vert s=0}+ \p_s\left(v\circ H_s^{-1}\right)(0,y_s,0)_{\vert s=0}=0\,.
\end{equation}
\end{lemme}

\proof A Taylor expansion gives 
$$
(v\circ H_s)(x_s,0,0)= (v\circ H_s)(0,0,0)+O(s^2)\,.
$$
That follows from the fact that  the velocity of $x_s$ is a vector which is contained in the kernel
of $dv$. Similarly, we have:
$$
(v\circ H_s^{-1})(0,y_s,0)= (v\circ H_s^{-1})(0,0,0)+ O(s^2)\,.
$$
Observe that $H_0(x,y,v)= (x,y,v)$. Thus, derivating the composed map  $H_s^{-1}\circ H_s= Id$
with respect to $s$ at $s=0$  yields:
$$\p_s H_s^{-1}(0,0,0)_{\vert s=0}+\p_s  H_s(0,0,0)_{\vert s=0}=0\,.
$$
Altogether, we get the desired formula.\bull

  \bigskip

\subsection{\sc Equators, signed hemispheres and 
 latitudes. \label{ssec:Equators}}${}$
 

We introduce some useful notations. Let $\bD^k$, { $ k\geq 1$,}
 be the closed Euclidean disc of dimension $k$ and radius 1 
  equipped with spherical coordinates $(r, \theta)\in [0,1]\x \bS^{k-1}$. A point $\theta\in \bS^{k-1}$
  will also be viewed as  unit vector $\theta\in T_0\bD^k$.
 
Suppose that we are given a preferred co-oriented hyperplane $\Delta\subset T_0\bD^k$. It determines  a \emph{preferred co-oriented equator} $E^\Delta\subset \bS^{k-1}$. 
The oriented normal to $\De$ determines two poles on the sphere: the \emph{North pole} $\nu_\De$ on the positive side 
of $\bD^k$ and the  \emph{South pole} $\si_\De$ on the negative side; 
and two open hemispheres of $\bS^{k-1}$ respectively noted $\cH^+(\bS^{k-1})$ and $\cH^-(\bS^{k-1})$. 


Any point $\theta\in\bS^{k-1}$ determines an angle with respect to the North pole $\nu_\De$. The cosinus of this angle defines a \emph{latitude} $\cos_\De:\bS^{k-1}\to [-1,1]$ defined by the scalar product
\begin{equation}\label{eq:Trigo}
\cos_\Delta(\theta):=\langle  \nu_\De,\theta\rangle.
\end{equation}

\begin{prop}\label{prop:Latitudes}
 Every $ X$ in $\S_g$ defines
 a preferred latitude on
 both the attaching sphere $\Si^-$ and the co-sphere $\Si^+$.
  \end{prop}

For this aim, we use 
 multispherical coordinates $(\phi,r,\psi)\in \bS^{i-1}\times  [0,1]\times \bS^{n-i-1}$ on 
 each level set of $\M_p$ (not well defined on the local stable/unstable manifolds). 
  We recall the map 
\begin{equation}\label{eq:Desc}
  Desc:\partial^+\M_p\smallsetminus \Si^+\to \partial^-\M_p\smallsetminus \Si^-
\end{equation}
obtained by descending the flow lines in $\M_p$. This map reads $Id$ in these coordinates.

The preferred latitude that we are going to define 
on $\Si^-$ and $\Si^+$ will be called respectively the \emph{$\phi$-latitude}
 and the \emph{$\psi$-latitude}.  We insist that these functions depend on $X\in \S_g$.
We denote them by 
\begin{equation}\label{eq:Latitudes}
 \cos_\phi^X:\Si^-\to [-1,1]\quad\text{and}\quad \cos_\psi^X:\Si^+\to [-1,1].
\end{equation} 
 When the vector field is clear from the context, these functions will just be denoted $ \cos_\phi$ and $ \cos_\psi$. 
 We shall decorate all the data 
related to $ \cos_\phi$ or $ \cos_\psi$ by using 
 the letter $\phi$ or $\psi$ respectively; 
 namely, the preferred hyperplane $\Delta^\phi$, the preferred equator $ E^\phi\subset \Si^-$, 
 the poles $ \nu_\phi$ and $\si_\phi $ in $\Si^-$, and so on.\\

\proof
Take any $X$ in $\S_g$ and denote by $\ell$ its  unique homoclinic 
orbit from $p$ to itself. The end point $a^+$ of $\underline\ell$
has coordinates $a^+=( -,0,\psi_0)$; as usual with polar coordinates, when the radius is 0 the spherical coordinate is not defined.
Let
\begin{equation}\label{eq:PsiProj}
\pi^{\psi_0}: \partial^+\M_p \to\{\psi=\psi_0\},\quad  \pi^{\psi_0}(\phi, r,\psi)= (\phi, r,\psi_0),
\end{equation}
  be the projection onto  the meridian $i$-disc.
 
  { Let $H: \{z=0\}\to \{z=1\}$ denote the holonomy diffeomorphism defined by the  vector field $X$
  in the tube $T$.
The image of $T\cap\Si^- $ through  $H$ 
 is a $(i-1)$-disc $D\subset \partial^+\M_p$.} Due to the transversality condition associated with $\ell$,
 this disc is a graph over 
 its projection $D_{\psi_0}:=\pi^{\psi_0}(D)$ { if the tube is small enough around $\underline\ell$.} 
Then, 
 \begin{equation}\label{eq:HypPhi}
\Delta^\phi:=T_{a^+}D_{\psi_0}\subset T_{a^+}\{\psi=\psi_0\}\text{ is the preferred hyperplane}.
\end{equation} 
As we noticed in  Remark \ref{rem:holonomie}, the vector $\p_v\in T_{a^+}\p^+\M_p$ is neither tangent to $\Si^+$ nor to $D$, which implies that
  \begin{equation} \label{positive_side}
  \parent{d\pi^{\psi_0}}_{a^+}(\partial_v) \text{ defines  { a co-orientation of }} \Delta^\phi  { \text{ in }
 T_{a^+}\{\psi=\psi_0\} .}
  \end{equation} 
This provides a preferred latitude on the $\phi$-sphere $\p\{\psi=\psi_0\}\cong  \bS^{i-1}\times\{1\}\times\{\psi_0\}$.
By the canonical isomorphism
  \begin{equation}\label{eq:IdentPhi}
\bS^{i-1}\times \{1\}\times\{\psi_0\} \cong \bS^{i-1}\times \{0\}\times\{-\} \cong \Si^-,
  \end{equation}
the preferred latitude on  $\p\sing{\psi=\psi_0}$ descends to some $\phi$-latitude 
which defines the announced $c^X_{\phi}:\Si^-\to [-1,1]$ in \eqref{eq:Latitudes}. The $\phi$-equator is the locus 
$\{c^X_\phi=0\}$, the North pole is $\nu_\phi=(c^X_\phi)^{-1}(1)$ etc.
 
For the $\psi$-latitude on the co-sphere $\Si^+$, we do the same construction by using
the reversed flow and its holonomy $H^{-1}$. More
precisely, 
take the image of $T\cap \Si^+$ through 
$H^{-1}$; it is a 
 $(n-i-1)$-disc $D'\subset \partial^-\M_p$ centered in $a^-$ whose spherical coordinates  are $ a^- = (\phi_0, 0, -) $.  Let 
 \begin{equation}\label{eq:PhiProj}
 \pi^{\phi_0}: \partial^-\M_p \to\{\phi=\phi_0\},\quad \pi^{\phi_0}(\phi, r,\psi)= (\phi_0, r,\psi),
 \end{equation}
 be  the projection onto the meridian $(n-i)$-disc
 and let  $D_{\phi_0}$ be the image $\pi^{\phi_0}(D')$.
Hence,
 \begin{equation}\label{eq:HypPsi}
\Delta^\psi:=T_{a^-}D_{\phi_0}\subset T_{a^-}\{\phi=\phi_0\} \text{ is the preferred hyperplane.}
\end{equation}

Moreover, $\p_v\in T_{a^-}\p^-\M_p$ is neither tangent to $\Si^-$ nor to $D'$, which implies 
that 
\begin{equation}\label{eq:PositSi+}
\parent{d\pi^{\phi_0}}_{a^-}(\partial_v)\text{ defines a co-orientation { of}  } \Delta^{\psi} { 
\text{ in } T_{a^-}\{\phi=\phi_0\}.}
\end{equation}
This yields a preferred latitude on the sphere $\p\sing{\phi=\phi_0}$, that can be carried to $\Si^+$
by means of
the canonical isomorphism
  \begin{equation}\label{eq:IdentPsi}
\sing{\phi_0}\x\{1\}\x \bS^{n-i-1}\cong \sing{-}\x\sing{0}\x \bS^{n-i-1} \cong \Si^+. 
  \end{equation}
This defines the announced $\psi$-latitude $ \cos_\psi^X$ in \eqref{eq:Latitudes}.
\bull

\subsection{\sc Holonomic factor and character function.}

By construction of the $\psi$-latitude and the $\psi$-latitude, we have the following splittings: 
\begin{equation}\label{eq:decompTgt}
\begin{cases}
T_{a^-}\{z=0\} = T_{a^-}\Si^-\oplus({\bR\nu_{\psi}}\oplus \Delta^\psi),\\
T_{a^+}\{z=1\} = (\Delta^\phi\oplus\bR{\nu_\phi})\oplus T_{a^+}\Si^+\, .
\end{cases}
\end{equation}

Given $X\in\S_g$, {recall  the homoclinic orbit $\ell$ whose homotopy class is $g\in \pi_1(M, p)$ and all
associated  objects that we 
introduced in Subsection \ref{ssec:tube}: the tube $T$, its coordinates $(x,y,v,z)$ 
and the holonomy diffeomorphism $H: \{z= 0\} \to  \{z= 1\}$. It reads $Id$
in the $(x,y,v)$-coordinates and $H(a^-)=a^+$.
We are free to choose the coordinates of the tube such that the 
unit tangent vector 
$\p_v^1:=\p_v\in T_{a^+}\{z=1\}$ verifies 
\begin{equation}\label{eq:ChoiceVAxis}
\p_v^1={\nu_{\phi}}.
\end{equation}

The \emph{linearized holonomy} $T_{a^+}H^{-1}$ 
maps $\p_v^1$ to $\p_v^0:=\p_v\in T_{a^-}\{z=0\}$.
By \eqref{eq:decompTgt}, the latter vector 
decomposes as 
\begin{equation}\label{eq:HolonDec}
\p_v^0=v_x+\bar{\eta}\,{\nu_\psi}+v_y,\text{ where }v_x\in T_{a^-}\Si^-, \, v_y\in \Delta^\psi,\, \bar{\eta}\in\bR.
\end{equation}
As we pointed in Remark \ref{rem:holonomie}, the only restriction
on the holonomy of $X$ along $\underline{\ell}$ is that $\bar{\eta}\neq 0$. 
Moreover, { according to \eqref{eq:PositSi+}, the vector  $\partial_v^0$ defines the positive side 
of the preferred hyperplane $\De^\psi$. As a consequence, }
 $\bar{\eta}$ must be positive.

\begin{defn}\label{def:HolFactor} 
The  holonomic factor associated with $X$ is the positive real number given by
\begin{equation}
\eta( X):=\,\frac{1}{\bar{\eta}}\,>\,0\, .
\end{equation}
\end{defn}
The following subsets of $\S_g$
 are respectively called the \emph{$\phi$-axis} and the \emph{$\psi$-axis} of $\S_g$:
 \begin{equation}\label{axis}
\S_g^{\phi}:=\ens{X\in\S_g}{a^+(X)\in E^\psi}\,\text{ and }\quad\S_g^{\psi}:=\ens{X\in\S_g}{a^-(X)\in E^\phi},
\end{equation}
that we also call the \emph{spherical axes}. Here, $a^+(X)$ and $a^-(X)$ 
stand for the respective extremities 
of the restricted orbit $\underline\ell$, where $\ell$ is the unique homoclinic orbit of $X$ in the homotopy class $g$.
Denote the intersection of the axes by:
\begin{equation} 
\S_g^{0,0}:=\S_g^\phi\cap\S_g^\psi
\end{equation}
which is empty when one axis is so. 

\begin{remarque} \label{empty-axis}When the Morse index of $\S_g$ is equal to 1, then the $\phi$-equator is empty but there 
are still signed poles. In that case, the $\phi$-latitude takes only the values $\{-1,+1\}$ and the $\psi$-axis 
is empty. When the Morse index of $\S_g$ equals $n-1$, then the $\psi$-equator and the $\phi$-axis are  empty  and the  $\psi$-latitude is valued in $\{-1,+1\}$. If $n>2$, these two events do not 
happen simultaneously. This is the reason for the dimension assumption in Theorem \ref{thm1} (2).

 \end{remarque}

We are now ready for defining the important notion of \emph{character function} togegther with its following 
ingredients.
In order to simplify notations, we introduce the {\it  extended $\phi$- and $\psi$-latitudes}
 to $\S_g$ by setting for every $X\in\S_g$: 
\begin{equation} \label{ext-latitude}
\om_\phi(X):= \cos_\phi^X(a^-(X)) \quad  \text{and}\quad \om_\psi(X):= \cos_\psi^X(a^+(X)).
\end{equation}

\begin{defn}\label{def:Character}
The \emph{character function} $\chi:\S_g\to\bR $  is defined by:
\begin{equation}\label{eq:CharVal}
\chi(X):=\eta( X)\,  \om_\psi(X)+ \om_\phi(X)
\end{equation}
Define  $\S_g^0$ (resp. $\S_g^+$, $\S_g^-$) as the locus where $\chi$ vanishes (resp. is positive, is negative). 
\end{defn}

By the very definition of the latitudes, 
it is clear that each axis intersects  $\S_g^0$ along $\S_g^{0,0}$, 
as Figure \ref{fig:Sg0etAxes} suggests (compare Figure \ref{fig:Sg0dsSg}).
\begin{figure}[h]
\includegraphics[width= 11cm, height=4cm]{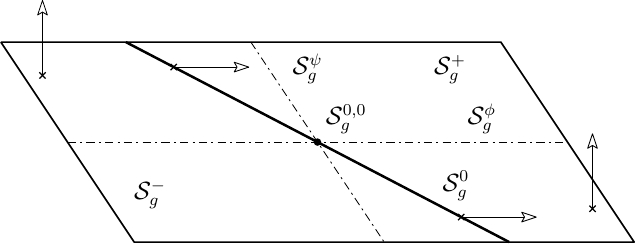}
\caption{The substratum $\S_g^{0,0}\subset \S_g^0$ as the intersection of the $\phi$-axis with the $\psi$-axis.} \label{fig:Sg0etAxes}
\end{figure}

Below, we start giving some information about $\S_g^{0,0}$ and $\S_g^0$  
from which Theorem \ref{thm1} will be completely proved.
\begin{prop}\label{item2} ${}$\\
\nd {\rm 1)} The axes $\S_g^\phi$ and $\S_g^\psi$  are 
$C^\infty$  submanifolds of codimension 1 in $\S_g$. Moreover, when  they are both non-empty their intersection
$\S_g^{0,0}$ is non-empty and transverse.
Hence, $\S_g^{0,0}$ is a $C^\infty$ submanifold of codimension 2 in $\S_g$. 

\nd {\rm 2)} If $n>2$ 
the  zero set $\S_g^0= \chi^{-1}(0)$ of the character function is
a \emph{non-empty} co-oriented $C^\infty$ submanifold of codimension 1 
 in \emph{each} connected component of $\S_g$. 
 \end{prop}
 
 \proof 
 Let $i=i(p)$ denote the Morse index of the zero $p\in Z(\al)$ where $g$ is based.\\
 1) The equation of the $\psi$-axis $\S_g^\psi$ in $\S_g$ reads 
 with the notations introduced in (\ref{ext-latitude}) and (\ref{axis}):
$$
  \om_\phi(X)= 0.
 $$
 If the index $i$ is equal to 1, by Remark \ref{empty-axis}, the $\psi$-axis is empty and there is nothing to prove. If not, 
 let $X_0\in \S_g^\psi$. We have to exhibit a germ 
 $(X_s)_{s\in \cO p(0)}$
 of path in $\S_g$ passing through  $X_0$ such that the $s$-derivative of 
 $  \om_\phi(X_s)$ at $s=0$ is non-zero. Let  $H_s$ be the local
 holonomy diffeomorphism
 of $X_s$ from a neighborhood of $\sing{z=0}$ in $\p^-\M_p$ to a neighborhood of $\sing{z=1}$ in $\p^+\M_p$. 
 Let $a^-=(\phi_0, 0,-)$ and $a^+=(-,0, \psi_0)$ be the end points of the restricted homoclinic 
 orbit $\underline\ell$ of $X_0$.
 We arrange that $H_s$ keeps the $\pi^{\psi_0}$-projection of $H_s(\Si^-)$ 
 into the { meridian $\{\psi=\psi_0\}\subset \p^+\M_p$ }
 independent of $s$. Thus, the equator $E^\phi$ is so and the $\phi$-latitude $\cos_\phi$
 does not depend on $s$. 
 Therefore, we are reduced to 
  control the $s$-derivative of $ \cos_\phi(a^-(X_s))$. 
 
 We recall that every germ of isotopy of the holonomy $H_0$ lifts to a deformation of $X_0$.
 Then,
 we are free to choose the holonomy so that $s\mapsto a^-(X_s)\in\Si^-$ crosses the non-empty equator $E^\phi$
 transversely at time $s=0$. 
 Thus, we are done.  
 For a similar reason, the equation $ \om_\psi(X)= 0$
of the { $\phi$-axis} $\S_g^\phi$  is regular. 

Let us  show the property of $\S_g^{0,0}$  when both axis are non-empty. In that case, the Morse index 
verifies  $1<i<n-1$, and it is available to have $\cos_\phi(a^-)=\cos_\psi(a^+)=0$, that is $\S_g^{0,0}\neq\emptyset$.
Let any $X_{0,0}\in \S_g^{0,0}$. 
We choose a 2-parameter family $X_{s,u}\in \S_g$ whose holonomy $H_{s,u}$
 satisfies the following conditions:
 \begin{enumerate}
 \item the equator $E^\phi$ is independent of $s$ when $u=0$ and $\p_s  \cos_\phi(a^-(X_{s,0}))>0$;
 \item  the equator $E^\psi$ is independent of $u$ when $s=0$ and $\p_u \cos_\psi(a^+(X_{0,u}))>0$:
 \item for every $(s,u)$ close to $(0,0)$, we have $a^-(X_{s,u})\in \Si^-$ and $a^+(X_{s,u})\in \Si^+$.
 \end{enumerate}
 Condition (3) guarantees that $X_{s,u}$ runs in $\S_g$. Thanks to (1) and (2), the evaluation map 
 $(s,u)\mapsto (a^-(X_{s,u}),  a^+(X_{s,u}))\in \Si^-\times\Si^+$ is transverse to the submanifold
 $E^\phi\times E^\psi$. This proves that the system of equations defining $\S_g^{0,0}$ near $X_{0,0}$,
 namely $ \left\{
  \om_\phi(X)=0, \ \om_\psi(X))=0
 \right\}$, is of  rank 2.\\
 
   \nd 2)  First, let us prove that the equation $\chi(X)=0$ has a solution in each connected component of $\S_g$.
   Let $X\in \S_g$ and let $a^+$ and $a^-$ be the corresponding end point of its restricted homoclinic 
   orbit $\underline\ell$. Any move of these points in their respective spherelifts to a deformation of $X$ in the space of adapted $\al$-gradients.
      If $\Si^+$ and $\Si^-$ are both connected, there is such a move until $a^-$ and $a^+$ lie in the 
    equators of their  respective sphere. Then, $X$ is deformed in $\S_g$ until it lies in $\S_g^{0,0}\subset \S_g^0$; this
    answer the question in this case.
    
    As $n>2$, one of the spheres $\Si^-$ and $\Si^+$ is not 0-dimensional, say $\Si^-\neq \bS^0$.
    Then, one can move $X$ in  $\S_g$ and modify the holonomic factor 
    by some homothety for making it  less than 1; secondly, knowing that $ \om_\psi(X)=\pm 1$,
    one moves $a^-$ in $\Si^-$ and changes $X$ accordingly, keeping the holonomic factor constant,
    up to reach the locus $\chi(X)=0$. So, $\S_g^0$ is visible in each connected component of $\S_g$.
    
    It remains to prove that the equation $\chi(X)=0$ is regular everywhere.  For every
    $X_0\in \S_g^0$ we have to exhibit 
   a germ of path $(X_s)_s$  in $\S_g$ passing through $X_0$ such that 
   $\p_s\chi(X_s)>0$.
  Let $\underline \ell$ be the restricted homoclinic orbit of $X_0$; let $a^+\in \Si^+$ and $a^-\in \Si ^-$ be its end points.
  First, we arrange that the equators $E^\phi $ and $E^\psi$ do not depend on $s$
  by requiring that the holonomy $H_s$ along the homoclinic orbit $\ell$ of $X_0$
   fulfils the following conditions:
   \begin{itemize}
   \item for every $s$, there is a homoclinic orbit $\ell_s$ (the end points of the restricted orbit $
   \underline\ell_s$ are noted $a^-(X_s)$ and $a^+(X_s)$);
   \item the $\pi^{\psi_0}$-projection of $H_s(\Si^-)$
   into the meridian $\{\psi=\psi_0\}$ is  constant; 
   \item the $\pi^{\phi_0}$-projection of 
   $(H_s)^{-1}(\Si^+)$ into the meridian $\{\phi=\phi_0\}$
    is constant. 
   \end{itemize}
   
   Now, there are two cases depending on whether $ \cos_\psi(a^+(X_0))$ is equal to 0 or not. 
   If $ \cos_\psi(a^+(X_0))$ is not 0, 
   the germ of $H_s$ at $a^-(X_0)$ is chosen to be a contraction: its
    center  is $a^+(X_0)$ 
   and its  factor  is $e^s$ (in the coordinate $(x,y,v)$ of the extremity $\{z=1\}$ of the tube $T$
   around $\underline\ell$). 
   Notice that such a contraction preserves the above requirements for the constancy of the equators. 
   Then, a calculation shows that the  holonomic factor
   is multiplied by the same factor, which implies that $\p_s\chi(X_s)>0$ since $a^\pm(X_s)$ is constant.
   
   Finally, we have to solve the case when  $ \cos_\psi(a^+(X_0))=0$. Here, we arrange the holonomy
   $H_s$ so that $\p_s \cos_\psi(a^+(X_s))>0$ and $\p_s \cos_\phi(a^-(X_s))=0$, which again 
   implies $\p_s\chi(X_s)>0$ since $\eta( X_s)>0$.  
  This finishes the proof of Proposition  \ref{item2}. \bull

  After Proposition \ref{item1} and Proposition \ref{item2} we have a complete proof 
 of Theorem \ref{thm1}.\bull

 \subsection{\sc Normalization of crossing path.}
 
 The normalization in question will be used  for proving Theorem \ref{thm:selfslideSimplif}
 and Theorem \ref{thm:Doubling}. The normalization is achieved by making some group act
 on $M$. 
  At the end of the subsection it will be proved that the stratification $(\S_g, \S_g^0, \S_g^{0,0})$
 is invariant under this action.
 
 In this subsection we use notations as $D_1(0), C_1(0)$ which will be used repeatedly 
   in Section \ref{section3} (see Notation \ref{not:C1(s)}). 
   Consider the image $H_0(\Si^-\cap\{z=0\})\subset \p^+\M_p$ 
     by the holonomy map of $X_0$
   along its homoclinic orbit $\ell$ in the homotopy 
   class $g\in \pi_1(M,p)$. Let us define 
$$D_1(0):= H_0(\Si^-)\cap\{z=1\}\ \text{ and }\  C_1(0):= Desc(D_1(0))$$
where $Desc$ is the descent map defined in (\ref{eq:Desc}). We recall $a^+=\ell\cap \Si^+$ 
(resp. $a^-= \ell\cap \Si^-$) whose spherical coordinate
is noted $\psi_0$ (resp $\phi_0$).

\begin{defn} \label{normalization} ${}$
A crossing path $(X_s)_{s\in \cO p(0)}$ of $\S_g$ is said to be normalized if it
fulfills the following requirements for every $s$, where $H_s$ denotes the perturbed holonomy of $X_s$.
\begin{enumerate}
\item  $D_1(0)$ has to be contained in
the preferred hyperplane $\De^\phi$ of the meridian $i$-disc\hfill\break $\{\psi=\psi_0\}$. 
\item Let $R_{a^-}$ denote the half ray $\{\phi= \frac \pi2 sign(\cos_\phi(a^-)),\, r\in [0,1], \psi=\psi_0\}$.
 The curve $J_0:=H_0^{-1}(R_{a^-})$ has to be contained in the meridian disc $\{\phi=\phi_0\}$ of $\p^-\M_p$.
\item The disc $D'(s):= H_s^{-1}(\Si^+)\cap\{z=0\}$ has to move in the same meridian disc $\{\phi=\phi_0\}$.

\end{enumerate}

\end{defn}

Only the first item will be used to prove Theorem \ref{thm:selfslideSimplif}; the two other items enter the proof of
Theorem  \ref{thm:Doubling}. The main tool for {\it normalization by conjugation} (see Proposition \ref{conjugation}) 
is given by the  next lemma about diffeomorphisms of $\M_p$.
Its proof by Taylor expansion is detailed in the Appendix to \cite{flX}.
  \begin{lemme} \label{C^1}
   Let $K$ be a $C^1$-diffeomorphism of $\partial^+\M_p$ of the form $(\phi, r,\psi)\mapsto
    (\phi, r, k(\phi,r,\psi))$ with
    $k(\phi,0,\psi)=\psi$. 
    Then, $K$ uniquely extends to $\M_p$ as a  $C^1$-diffeomorphism
    which  is the identity on both stable and unstable local manifolds 
    and which  keeps 
    the standard vector field $X_i$ invariant. Moreover, the extension ${\overline K}$
    is $C^1$-tangent to $Id$ along the attaching sphere $\Si^-$.
   \end{lemme}
   
   It is worth noting that the extension cannot be $C^2$ in general, even if $K$ is $C^\infty$. 
   This lemma can be also used by interchanging the roles of $\partial ^+\M_p$ and $\partial ^-\M_p$
   and simultaneously the roles of $\phi$ and $\psi$.

   \begin{prop}\label{conjugation} Given a positive crossing path $(X_s)_s$ of the stratum $\S_g$,
   there exists a $C^1$-diffeomorphism $\overline K$ of $M$,  isotopic to $Id_M$ among the 
   $C^1$-diffeomorphisms keeping $\al$ and $\M_p$ invariant, 
   such that the crossing path $(\overline K_* X_s)_s$ carried by $\overline K$ is normalized.
    Moreover, $\overline K$ may be chosen so that it preserves $\ell$ pointwise.
   \end{prop}
   
   Notice that the vector field $(\overline K_* X_s)$ might only be $C^0$. But it is integrable
   and the associate foliation is  transversely $C^1$;
   its holonomy is changed by $C^1$-conjugation.\\ 
   
   \proof 
   \nd 1) We first look for a diffeomorphism $\overline K$ of $M$ carrying $(X_s)_s$ to a crossing path which fulfils
   the first  item of Definition \ref{normalization}.
     If the tube $T$ around $\ell$ is small enough $D_1(0)$ is nowhere tangent to the fibres 
   of the projection  $\pi^{\psi_0}$ 
   to the meridian disc $\{\psi=\psi_0\}$. As a consequence,
   its projected disc $D_{\psi_0}$ is smooth and there exists a smooth map $\bar k:D_{\psi_0}\to \Si^+$
   such that  $D_1(0)$ reads  
   $$D_1(0)= \{(\phi,r,\bar k(\phi,r))\mid (\phi,r)\in D_{\psi_0}\}.$$
   Since its source is contractible,  $\bar k$ is homotopic to the constant map valued in $\psi_0$.
   By  isotopy extension preserving the fibres of $\pi^{\psi_0}$, there exists some 
   diffeomorphism  $K_1$ of $\M_p$ of the form assumed in
    Lemma \ref{C^1} which maps the given $D_1(0)$   to $D_{\psi_0}$. Therefore,
   this $K_1$ extends to $\M_p$ 
   Since $K_1$ is isotopic to $Id$ through diffeomorphisms of the same type, its extension
   $\overline K_1$ to $\M_p$ also extends to $M$ with the same name. Moreover, 
  the isotopy of $\overline K_1$ to $Id_M$ is supported in a neighborhood of $\M_p$ 
  and preserves each level  set of a local primitive of $\al$. Since $\Si^\pm$ are kept fixed by 
  $\overline K_1$, it is easy to get that $\ell$ is fixed by $\overline K_1$.
  
  After having carried $X_0$ by this $\overline K_1$, we are reduced to the case where 
  $D_1(0)$ is contained in the meridian disc $\{\psi=\psi_0\}$. Decreasing the radius
  of the tube $T$ if necessary, the tangent plane 
  $T_mD_1(0)$ is almost orthogonal to the pole axis directed by $\nu_\phi$
  at each $m\in D_1(0)$. This implies that, for every $r\in (0,1)$, the disc $D_1(0)$ is transverse
  to the $(i-1)$-sphere of radius $r$ in $\{\psi=\psi_0\}$.
   
   The image $C_1(0)$   of $D_1(0)$ by $Desc$
       is diffeomorphic to  $S^{i-2}\x(0,1] $ 
  and contained in the spherical annulus $\A_{\psi_0}:=\{(\phi, r, \psi_0)\mid \phi\in \Si^-, r\in (0,1]\}$. By tangency of $D_1(0)$ with the preferred hyperplane $\De^\phi$,
  the 
 {\it end} of $C_1(0)$ when $r\to 0$   compactifies as the 
   $\phi$-equator  $E^\phi\subset \Si^-$. 
  Moreover, $C_1(0)$ is transverse (inside $\A_{\psi_0}$) to the sphere $\Si^-\times\{(r,\psi_0)\}$ for every
   $r\in(0,1)$, since the corresponding assertion holds in $\p^+\M_p$.
  Thus, there is an annulus $C_{\rm eq}\subset E^\phi\times[0,1]\times\{\psi_0\} $
  such that $C_1(0)$ reads as the graph of some map $\bar\kappa: C_{\rm eq}\to \Si^-$ valued in the 
  complement of the poles. Then, $\bar\kappa$ is homotopic to the map $(\phi,r)\mapsto \phi$
  from $C_{\rm eq}$ to the equator $E^\phi$ of $\Si^-$. 
  
  By isotopy extension preserving  each sphere $\Si^-\times\{(r,\psi)\}$,
   we have some diffeomorphism  $\overline K_2$ of $\p^-\M_p$ of the form 
   $(\phi,r, \psi)\mapsto (\kappa(\phi,r,\psi),r,\psi)$ which
    pushes $C_1(0)$ to its {\it flat} position $C_{\rm eq}$ and satisfies $\kappa(\phi,0,\psi)=\phi$.
   By applying Lemma \ref{C^1} ``up side down'', $K_2$ extends to $M$ preserving $\M_p$ with its 
   standard gradient. On the upper boundary of the Morse model,
  this means that $\overline K_2$ pushes $D_1(0)$ to an $(i-1)$-disc in the hyperplane $\De^\phi$
  by some diffeomorphism tangent to $Id$ in $a_+$. As for $\overline K_1$, this $\overline K_2$
  may be chosen so that $\ell$ is fixed pointwise. The composed diffeomorphism 
  $\overline K_2\circ\overline K_1$ is as desired.\\

  \nd 2) We now prove the last two items. It consists just in an easy addition to what we have done above. 
  The diffeomorphism $\overline K_1$ keeps $\ell$ fixed pointwise. Let $\pi^{\phi_0}$ be the projection
  of $\p^-\M_p\cong \left(\Si^-\times \{\phi=\phi_0\}\right)$ onto its second factor. The tangent space in $a^-$
  to $\Si^-$, $D'(0)$ and $J_0$ are independent and their span is egal to $T_{a^-}\p^-M_p$. Then, after shrinking the 
  tube $T$ if necessary $\pi^{\phi_0}$ embeds the union $\left(\cup_sD'(s)\right)\cup J_0$ into $ \{\phi=\phi_0\}$;
  notice that due to the non-vanishing of the velocity with respect to $s$ the union $\left(\cup_sD'(s)\right)$
  is an $(n-i)$-disc transverse to the fibres of $\pi^{\phi_0}$. 
  
  As $a^-$ is far from the equator $E^\phi$, it is easy to find a common diffeomorphism $K_2$ of $\p^-\M_p$
  preserving the coordinates $(r,\psi)$  such that it maps $C_1(0)$ to the equatorial annulus $C_{\rm eq}$---the job 
  required for Item (1)---and 
  simultaneously $\left(\cup_sD'(s)\right)\cup J_0$ onto its  $\pi^{\phi_0}$-image in the meridian $\{\phi=\phi_0\}$ by 
  diffeomorphism.
 This complete the proof. 
   \bull

\begin{remarque} Proposition \ref{conjugation} holds true for a  finite dimensional family. For instance,
if we are given a two-dimensional germ $\left(X_{s,t}\right)_{s,t}$ adapted to the pair $(\S_g,\S_g^0)$
 in $X_{0,0}$---in the sense of Definition \ref{decompositionS^0}---then each crossing path $\ga_t:=(s\mapsto X_{s,t})$  of $\S_g$ has a normalization 
 by some $C^1$-diffeomorphism $\overline K_t$ depending continuously on $t$ in the $C^1$-topology.
 \end{remarque}

\begin{notation}\label{cG}
Let $\cG^\pm$ be the groups of diffeomorphisms of $M$ isotopic to $Id_M$ which  fixe the homoclinic orbit
$\ell$ pointwise, preserve the  closed one-form $\al$ and its standard gradient in $\M_p$,
and have  the following form:
\begin{itemize}
\item the restriction of every element  in $ \cG^+$  to $\p^+\M_p$ read $\bparent{\phi,r,\psi)\mapsto (\phi, r, k^+(\phi,r, \psi)}$
with $k^+(\phi, 0,\psi)=\psi$;
\item the restriction of every element in $ \cG^-$ to $\p^-\M_p$ 
reads $(\phi,r,\psi)\mapsto \bparent{k^-(\phi, r,\psi),r, \psi}$
with $k^-(\phi, 0,\psi)=\phi$.
\end{itemize} 
\end{notation}

\begin{prop} \label{invariance}
The action of the  groups $\cG^+$ and $\cG^-$  on the space of adapted $\al$-gradients
preserves the strata $\S_g$, $\S_g^0$ and $\S_g^{0,0}$.

\end{prop}

  \proof   We do it for $\cG^+$. Let $G\in \cG^+$ and $X_0\in \S_g$. Since
$G$ fixes the homoclinic orbit $\ell$ pointwise, the carried 
vector field $G_*(X_0)$
has the same homoclinic orbit. According to the form of the restriction of $G$ to the 
upper boundary of $\M_p$, the projection of $D_1(0)$ to the meridian disc  is unchanged.
Therefore, the $\phi$-equator is preserved. Looking in the lower boundary, one derives that the
$\phi$-latitude of $a^-(X_0)= \Si^-\cap\ell$ is preserved.


Consider the disc $D'_1(0):= H^{-1}_{X_0}(\Si^+)$. Recall from Lemma \ref{C^1} that $G\vert_{\p^-\M_p}$ is tangent to $Id$ at every point of $\Sigma^-$. Therefore, the tangent space $T_{a^-}D'_1(0)$ remains invariant by $G$. It follows that the $\psi$-equator is not changed, and hence,
 the $\psi$-latitude of $a^+$ is preserved. Thus we have the invariance of the spherical 
 axes and of their intersection $\S_g^{0,0}$.
 
 It remains to show that the character function is invariant. We already have seen the invariance of the 
 latitudes. The last term to control is the holonomic factor $\eta(X_0)$---resp. 
 $\eta(G_*(X_0))$---defined in (\ref{eq:HolonDec}). This factor remains unchanged by the action of $G$ thanks to the invariance of:
 \begin{itemize}
\item $\p_v^1$ by invariance of the $\phi$-latitude (see Equation (\ref{eq:ChoiceVAxis}),
\item $\p_v^0$ since 
 $DG_{a^-}=Id$,
\item the framing in which $\p_v^0$ decomposes (this framing is preserved by invariance of the $\psi$-latitude). \end{itemize}
  \bull

 \section{ Change in the Morse-Novikov complex}\label{section3}
  \subsection{\sc A groupoid approach.} \label{ssec:Groupoid}

 A groupoid $\sG$ is a {\it small category} where every arrow is invertible. The set of objects in $\sG$ is noted $\sG^0$
 and the set of arrows (or morphisms) is noted $\sG^1$. Given two objects $p,q\in \sG^0$, the set of arrows
 from $p$ to $q$ is noted $\Hom(p,q)$. 
 The  \emph{identity} arrow at $p$ is noted $1_p\in\Hom(p,p)$. The map $p\mapsto 1_p$ embeds $\sG^0$
  into $\sG^1$; and hence, $\sG$ may be identified with its set of arrows endowed with its subset of identity arrows.
  The maps \emph{source} and \emph{target}, $s,t:\sG^1\to\sG^0$ are defined by $s(g)= p$ and $t(g)= q$ for every morphism 
 $g\in \Hom(p,q)$.
 
\begin{remarque}\label{rem:GroupoidRing}
 We denote by $\Z[[\sG]]$ the set of formal series of the arrows of $\sG$. An element $\la\in\Z[[\sG]]$
is usually written as $\la=\sum_{g\in\sG}n_g(\la)g$, where $ n_g(\la)\in\Z$. Define the \emph{support} of $\la$ as the set $\supp(\la):=\ens{g\in\sG}{n_g(\la)\neq 0}$. Consider the set
\begin{equation}
\Z[\sG]:=\ens{\la\in\Z[[\sG]]}{\supp(\la)\text{ is finite}}.
\end{equation}
Given two arrows $g,h$---seen as elements of $\Z[\sG]$--- the product $gh\in\Z[\sG]$ is defined by their composition in $\sG^1$ when $t(g)=s(h)$ and by $0$ otherwise. Extending the previous rule distributively with respect to the sum, we obtain a ring structure for $\Z[\sG]$. Moreover, when $\sG^0$ is finite, the element $1 := \sum_{p\in\sG^0}1_p\in\Z[\sG]$ gives an \emph{identity element} for this product. We call $\Z[\sG]$ the \emph{groupoid ring} associated with $\sG$. The next definition is classical.
\end{remarque}

\begin{defn}\label{fund_gr}
The fundamental groupoid $\Pi$ of the manifold $M$ is defined as follows: its objects are the points of 
$M$ and if $(p,q) $ is a pair of points $\Hom(p,q)$ is the set of  homotopy classes
of paths from $p$ to $q$. If $\ga$ is a such a path, its homotopy class $[\ga]$ will be called the 
$\Pi$-value of $\ga$.
\end{defn}

 The closed 1-form $\al$ (whose cohomology class is noted $u$) defines a groupoid morphism
\begin{equation}\label{u-alpha}
 u_\al:\Pi\to\R,  
\quad\quad g\mapsto \int_{\gamma}\alpha,
\end{equation}
 where $g$ is the $\Pi$-value
 of a path $\gamma$ in $M$ and $\R$
  is seen as a groupoid with a single object $0$.
   The restriction of any such $u_\al$  to the fundamental group $\pi_1(M,p)$ clearly coincides 
   with the group morphism associated with $u$. \\

We denote by $\Pi_\al$ the full subcategory of $\Pi$ whose set of objects 
 is the set $Z(\al)$, the zero set of $\al$. By Remark \ref{rem:GroupoidRing}, when $\al$ is Morse and $Z(\al)$ is non-empty, we may consider the groupoid ring $\Z[\Pi_\al]$.\\

A formal series $\la\in\Z[[\Pi_\al]]$ fulfills the \emph{Novikov Condition} if
\begin{equation}\label{eq:NovikovCond}
\text{for every } L\in\R, \text{ the set } \supp(\la)\cap u_\al^{-1}(L,+\infty) \text{ is finite}.
\end{equation}

Denote by $\La_u\subset \Z[[\Pi_\al]] $ the subset of formal series satisfying the Novikov Condition. It can be easily checked that the product rule given in Remark \ref{rem:GroupoidRing} also gives a ring structure to $\La_u$, having the same identity element as $\Z[\Pi_\al]$. We call $\La_u$ the \emph{Novikov ring associated with $\al$}.

\begin{ex}
Let $p\in  Z(\al)$ and $ g\in \pi_1(M, p)$  with $u(g)<0$ (for instance, 
 the $\Pi$-value of a homoclinic orbit of some $\al$-gradient).
The following formal series are elements of the Novikov ring:
\begin{equation}
\sum_{j=1}^\infty  g^j\,\quad\text{and}\quad 
\sum_{j=1}^\infty (-1)^j g^j\
\end{equation}

Indeed, the Novikov Condition  \eqref{eq:NovikovCond} is fulfilled
since $u_\al(g^j)=j.u_\al(g)$ which goes to $-\infty $
as $j\to +\infty$. Thus, these two series belong to the Novikov ring $\La_u$.
 In particular $1+g+g^2+\ldots$ is a 
 unit whose inverse is $1-g$. \\
 \end{ex}

\subsection{\sc The Morse-Novikov complex.} 
Let $X$ be an $\al$-gradient which is assumed $KS$. An orientation is arbitrarily chosen 
on the unstable manifold $W^u(p,X)$ for each zero $p\in Z(\al)$.
We are going to define a chain complex $\bparent{C_*(\al),\p^X}$ of $\La_u$-modules; 
it is graded by the integers $i\in\{0, 1, \ldots, n=\dim M\}$.
This complex will be called the {\it Morse-Novikov} complex associated with the pair $(\al, X)$.

For each degree $i$, the module $C_i(\al)$ is the left $\La_{\al}$-module freely generated by $Z_i(\al)$, the finite set 
of zeroes of $\al$ of Morse index $i$.
The $\La_u$-morphism $\partial_* ^X: C_*(\al)\to C_{*-1}(\al)$  must have the form following form 
on each generator of $C_i(\al)$: 
 \begin{equation}\label{eq:NovDiff}
\p^X_*(p)=\sum_{q\in Z_{*-1}(\al)}n(p,q)^X q,
\end{equation}
where the coefficient of $q$ has to be an element of $\La_u$ (called the 
 \emph{incidence} of $p$ to $q$). This coefficient
 $n(p,q)^X$ is the algebraic count 
  which we are going to define.
  
   Let $\Orb^X(p,q)$ denote the set of connecting orbits of $X$
  from $p$ to $q$. 
 First, we define the sign of a connecting orbit $\ell\in\Orb^X(p,q)$. Given a point $x\in\ell$,
the sign $\ep_\ell$ is defined by the following equation: 
\begin{equation}
\ep_\ell\, X(x)\wedge \text{co-or}\parent{W^s(q,X)}= \text{or} \parent{W^u(p,X)}.
\end{equation}
This definition is clearly indepedent of $x\in\ell$.}
\begin{defn} 
Assume the $\al$-gradient $X$ is $KS$.
The \emph{Morse-Novikov incidence} associated with the data $(p,q,X)$, $p\in Z_i(\al)$,
$q\in Z_{i-1}(\al)$,
 is defined by:
\begin{equation}
n(p,q)^X:=\sum_{\ell\in\Orb^X(p,q)}\ep_\ell \,g_\ell 
\end{equation}
where  $g_{\ell}$ denotes the $\Pi$-value of the connecting orbit $\ell$.
\end{defn}

 By Proposition \ref{deep-appendix}, this coefficient fulfills the Novikov Condition (\ref{eq:NovikovCond}).
 So, it is an element of $\La_u$. Moreover, 
the map $\p^X$ as in \eqref{eq:NovDiff} is indeed a differential; this can be found in \cite{latour}. The resulting $\bparent{C_*(\al),\p^X}$ is known as the Morse-Novikov complex (see \cite{novikov}, \cite{sikorav}).\\

We denote by $\Orb^X_L(p,q)$ the set of connecting orbits from $p$ to $q$ whose $\al$-length $ \cL(\ell)$ is 
less than a fixed $L>0$. Since these orbits verify the inequality $u_{\al}(g_\ell)>-L$, we are led to define a 
$L$-\emph{truncation} map $\cT_L: \La_u\to\Z[\Pi_\al]$ by:
\begin{equation}
\cT_L(\la):=\sum_{u_{\al}(g)>-L}n_g(\la)g.
\end{equation}

Two elements $\la,\mu\in\La_{\al}$ are said to be  equal modulo $L$ 
if $\cT_L(\la-\mu)=0$.\\

Finally, the \emph{$L$-incidence} is defined as follows: \begin{equation}
n(p,q)^X_L:=\sum_{\ell\in\Orb^X_L(p,q)}\ep_\ell \,g_\ell\in\Z[\Pi_\al].
\end{equation}
Of course
we have $\cT_L\bparent{n(p,q)^X}=n(p,q)^X_L$.

   \subsection{\sc Effect of homoclinic bifurcation on the incidence.}  
   
 Consider a generic one-parameter family of addapted $\al$-gradients  $(X_s)_s$ 
 such that $ X_0\in\S_g$. By definition, 
$X_0$ has a unique homoclinic orbit $\ell$ connecting  $p\in Z(\al)$ to itself 
whose   $\Pi$-value is $g$. Denote the index of $p$  by $i$.

The next definition specifies some genericity conditions that will be needed to prove the theorem below. The rest
of this section is devoted to its proof and consequences.
 
 \begin{defn} \label{almost}Let $L>0$.

\nd{\rm 1)} The $\al$-gradient  $X$  is said to be \emph{Kupka-Smale up to} $L$ if, for every pair of 
zeroes $p,q\in Z(\al)$ and every $X$-orbit $\ell$
from $p$ to $q$ with $-\int_\ell \al<L$, the unstable and stable manifolds, $W^u(p,X)$
and $W^s(q,X)$, are transverse along $\ell$. The subset of $\F_\a$ formed with such $\al$-gradients 
is noted $\mathcal K S_L$.
\smallskip

\nd{\rm 2)} An $\al$-gradient $X\in \S_g$ 
 is said to be 
\emph{almost Kupka-Smale up to }$L$,
if the preceding transversality condition is fulfilled except for the unique homoclinic orbit
whose $\Pi$-value is $g$. The subset of $\S_g$ formed with such elements is noted $\S_{g,L}$.\smallskip

\nd{\rm 3)} The $\al$-gradient $X\in \S_g$
 is said to be \emph{almost Kupka-Smale} if it is KS up to L for every $L>0$. These gradients
are the elements of  $\S_{g,\infty}: = \bigcap_L\S_{g,L}$.
 \end{defn}

\begin{prop}${}$\label{g-residual}

\nd{\rm 1)} The subspace $\S_{g,L}$ is open and dense in $\S_g$. Moreover, there exists an open set $\mathcal W_L$
in $\F_\al$  such that $\mathcal W_L\smallsetminus \S_g$ is contained in $\mathcal K S_L$.

\nd{\rm 2)} The subspace $\S_{g,\infty}$ is residual in $\S_g$.
\end{prop}

\proof 1) One checks that the constraint to have 
a unique homoclinic orbit with a given $\Pi$-value does not prevent us from arguing as Peixoto \cite{peixoto}.

\nd 2) This item is a little more subtle since $\S_g$ is not a complete metric space. But it is separable. Then, it is 
sufficient to prove that, for any closed ball $B$ in $\S_g$, the intersection $\S_{g,\infty}\cap B$ is residual
in $\S_g\cap B$. And now, we are allowed to follow Peixoto word for word. More details are left to the reader. \bull

\begin{notation}{\rm It is easily seen that $\S_{g,L}$ is open in $\S_g$; and, if $X\in \S_{g,L}$
 there is an arbitrarily small neighborhood $U$ of $X$ in $\F_\al$
such that $U\smallsetminus \S_g$ is made of two connected components in $\mathcal K S_L$.  
In particular, if $(X_s)_s$ is a path which intersects $\S_g$ transversely at $ X_0\in \S_{g,L}$ the gradient $X_s$
 is  $KS$ up to $L$ for every $s\neq 0$ close enough  to 0.
 
 Therefore, for every $q\in Z_{i-1}(\al)$
 the $L$-incidence $n(p,q)^{X_s}_L$ is well defined and independent of $s$
when $s\in \cO_+$ (resp. $s\in \cO_-$); 
it is  denoted by $n(p,q)^\pm_L$ respectively. Here, the symbol $\cO_-$ stands for an open interval $(-\epsilon,0)$ 
whose size is not specified and which is as small 
as needed by the statement; and similarly for $\mathcal O_+$.
}
\end{notation}

\begin{thm}\label{thm:selfslideAlg} 
Let $(X_s)_{s}$ be a path of adapted $\al$-gradients which intersects 
 $\S_g$ transversely at $ X_0$ and let $L> - u_\al(g)$.
Assume $X_0$ is almost Kupka-Smale up to $L$, that is $X_0\in \S_{g,L}$.
Then we have the following.
\smallskip

\nd  When $(X_s)_s$ intersects the stratum $\S_g$ positively, 
 the next relations hold in $\La_u$: 
\begin{enumerate}
\item if $ X_0\in\S_g^+$, then $n(p,q)^+_L=\bparent{1+g+g^2+g^3+\ldots}
\cdot n(p,q)^-_L$\quad{\rm (mod. $L$)},
\item if $X_0\in\S_g^-$, then $n(p,q)^+_L=\bparent{1+g}\cdot n(p,q)^-_L$\quad{\rm (mod. $L$)}.
\end{enumerate}
\smallskip

\nd
When $(X_s)_s$ intersects the stratum $\S_g$ negatively, we have:
\begin{enumerate}
\item[(1')] if $ X_0\in\S_g^+$, then $n(p,q)^+_L=\bparent{1-g}\cdot n(p,q)^-_L$\quad{\rm (mod. $L$)},
\item[(2')] if $ X_0\in\S_g^-$, then $n(p,q)^+_L=\bparent{1-g+g^2-g^3+\ldots}\cdot n(p,q)^-_L$\quad{\rm (mod. $L$)}. 
\end{enumerate}

\end{thm} 

It is worth noticing the reason for the truncation: in general, the bifurcation at $s=0$ is not isolated among the bifurcation times of the path
$(X_s)_s$. When it is isolated, the truncation is not needed any more; this will be the case in  
\cite{l-m}.\\

\nd{\bf Proof of $(1)\iff(1')\iff (2)\iff (2')$.} The first equivalence is obvious 
since $1- g$ is the inverse of $ 1+\Si_{j=1}^\infty g^j$; and similarly for the last equivalence. 

Let us show the middle equivalence. It is obtained by changing the vector $\partial_v$ 
into its  opposite in the coordinates of the tube around the homoclinic orbit of $X_0$. This amounts to put a sign in  Formula (\ref{orientation1}). The latter change has three effects: 
\begin{enumerate} 
\item[i)] It reverses the co-orientation of $\S_g$. Hence, positive and negative crossings are exchanged.
\item[ii)] The  character is changed into its opposite since the $\phi$- and $\psi$-latitudes are.
Thus, both sides of $\S^0_g$ are exchanged. 
\item[iii)] The homoclinic orbit becomes negative in the following sense: if the $\phi$-sphere
is seen 
as the boundary of the meridian disc at 
$a^+$,  the new positive hemisphere  projects to the preferred
hyperplane $\De^\phi$ (whose orientation is unchanged) by reversing the orientation. This implies that in the algebraic count 
of connecting orbits from $p$ to $q$ (where $i(p)=i(q)+1$) the coefficient 
$g$ has into be changed to $-g$ (see
the orientation claim in Lemma \ref{lem:C1(s)}).
\end{enumerate}
We are left to prove the Theorem \ref{thm:selfslideAlg}  in case (1). This will be done 
in Subsection \ref{continued}. 
\bull


According to Proposition \ref{invariance}, the statement of Theorem \ref{thm:selfslideAlg}
 is invariant by the groups $\cG^\pm$ 
introduced in Notation \ref{cG}. After Proposition \ref{conjugation}, it is sufficient to consider the case
where the crossing path in question is {\it normalized} in the sense of Definition \ref{normalization}. 
This assumption is done in what follows.
We need some more notations and lemmas.
The setting of Theorem \ref{thm:selfslideAlg} is still assumed.

\begin{notation}{\rm ${}$\label{not:C1(s)}

\nd 1) Recall from Subsection \ref{perturbed} that, $H_s$ denotes the perturbed holonomy diffeomorphism
 along the homoclinic orbit $\ell$. For $s\in \cO p(0)$, it maps  $\mathcal Op(\{z=0\}) \subset \partial^- \M_p$ to an 
 open set of $\partial^+ \M_p$ containing $\{z=1\}$.
\smallskip

\nd 2) For $s\in \mathcal Op(0)$, 
let $D_1(s)$ denote the image  $H_s(\Si^-)\cap\{z=1\}$. Consider its projection  $\pi^{\psi_0}(D_1(s))$ 
onto the meridian  disc $\{\psi=\psi_0\}$ and define 
\begin{equation}\label{a+(s)}
a^+(s):=  \pi^{\psi_0}(D_1(s))\cap\R\partial_v^1. 
\end{equation}
The {\it crossing velocity} of the crossing path is 
\begin{equation}\label{initial_velocity}
V_1 := \frac {d }{ds} a^+(s)_{\vert s=0}. 
\end{equation}
After reparametrization, we may assume $V_1=1$.
\smallskip

\nd 3)  For $s\in \mathcal O_\pm$, by definition of a crossing path $D_1(s)$ avoids $\Si^+$. Therefore we are allowed  to define $C_1(s):= Desc(D_1(s))\subset \partial^-\M_p$. It is still
an $(i-1)$-disc.
}
\end{notation}

\begin{lemme}\label{lem:C1(s)}
Recall the natural projection $\pi^-: \partial^- \M_p\to \Si^-$.
Let $K$ be any compact disc in the open hemisphere $\mathcal H^-(\Si^-)$ {\rm (}as in {\rm  Subsection \ref{ssec:Equators})}. 
Then, for $s\in \mathcal O_-$, the disc $C_1(s)\cap (\pi^-)^{-1}(K)$ is a graph over $K$
of a section of $\pi^-$. 
This section goes to the zero-section $0_K$ of $\pi^- $ over $K$
in the $C^1$-topology as 
 $s$ goes to $0$. 
Moreover, $\pi^-: C_1(s)\to \Si^-$ is orientation reversing. 

A similar statement holds when $K\subset \mathcal H^+(\Si^-)$ and $s\in\mathcal O_+$,
except that $\pi^-: C_1(s)\to \Si^-$ is orientation preserving in that case.
\end{lemme}

\proof 
The statement about orientation is clear after the claim 
about the $C^1$-convergence.
Consider the case $s\in\mathcal O_-$, the other case being similar. 
 Recall the normalization assumption: the disc $D_1(0)$ is contained in the meridian disc $\{\psi=\psi_0\}$.
 Recall the projection $\pi^{\psi_0}$ of $\p^+\M_p$ to the meridian disc.
The normalization implies that the projected discs $\pi^{\psi_0}(D_1(s))$ tend to  $D_1(0)$
in the $C^1$-topology.

  Recall the identification $\p(\{\psi=\psi_0\})\cong \Si^-$ of \eqref{eq:IdentPhi}
  and think of $K$ as a compact subset of the South hemisphere in the boundary of the meridian disc
  $\{\psi=\psi_0\}$. 
  For every such $K$, 
   the following property holds:
  \begin{itemize}
  \item[] {\it For every $s$ close enough to 0 and for every $\phi\in K$, 
  the disc $\pi^{\psi_0}(D_1(s))$
  intersects only in one point and transversely the ray directed by $\phi$ in $\{\psi=\psi_0\}$. }
  \end{itemize}
   This point is denoted by $m_s(\phi)$; it is the image of some 
   $\tilde m_s(\phi)\in D_1(s)$ through $\pi^{\psi_0}$. 
   We have $m_0(\phi)=\tilde m_0(\phi)= a^+$, but when $s\neq 0$, the point 
   $\tilde m_s(\phi)$ has well-defined multi-spherical coordinates $(\phi, r_s(\phi),\psi_s(\phi))$ where
   $r_s$ and $\psi_s$ depend smoothly on $s$.

 Going back to $\p^-\M_p$ by the map $Desc$, we see that  $C_1(s)$ is  the image of a section 
 of the projection $\pi^-$ over $K$. 
 When $s\in \cO_-$ goes to 0, then $D_1(s)\cap \{(\phi,r,\psi)\in \p^+\M_p
 \mid \phi\in K\} $ 
  goes to $a^+$  
in the metric sense. In particular, $\max_{\phi\in K} \{r\,\vert\, (\phi,r,\psi)\in D_1(s)\}$ 
 goes to 0.
Therefore, $C_1(s)\cap (\pi^-)^{-1}(K)$ goes to $0_K$ in the $C^0$-topology
when $s$ goes to  $0$ negatively. 

For the $C^1$-convergence, we use  that $K$ is far from the 
$\phi$-equator of $\Si^-$.
Therefore, the angle in the meridian disc  $\{\psi=\psi_0\}$ between the ray directed by  $\phi$ 
and the tangent plane to $ \pi^{\psi_0}(D_1(s))$ at  $ m_s(\phi)$ 
is bounded from below.  Including the fact that $s\to 0_-$ implies $r\to 0$, it follows 
that the smooth functions $r_s(\phi)$ and $\psi_s(\phi)$ satisfy
$$\left\{\begin{array}{c}
\vert dr\vert= O_s(r) \vert d\phi\vert,\\
\vert d\psi\vert =O_s(r) \vert d\phi\vert,
\end{array}
\right.
$$
where $O_s(r)$ stands for a quantity which is uniformly bounded by a positive multiple of  $r$ when $s$ goes to 0.
This yields the claimed 
$C^1$-convergence 
of the part of $C_1(s)$ over $K\subset \mathcal H^-(\Si^-)$ to $0_K = K$.
\bull

\subsection{\sc Geometric interpretation of the character function.}\label{interpretation}
We still consider a germ of {\it normalized} positive crossing path $\left(X_s\right)_s$.
Let $D_1'(0) $ be the connected component $W^s(p, X_0)\cap \{z=0\}$ which contains $a^-$. 
This is an $(n-i-1)$-disc which is the image of $\Si^+$  
by the inverse holonomy diffeomorphism $H_0^{-1}$ along $\ell$.
For every $s\in \mathcal Op(0)$, consider now $D_1'(s):= H_s^{-1}(\Si^+)\cap \{ z=0\}$.

Recall from Subsection \ref{ssec:tube}  
 that  $\Si^-\cap  \{z=0\}$ is identified 
with the $x$-axis whereas $D_1'(0)$ is identified with the $y$-axis. 
 Let also $p_v:\sing{z=0}\to\sing{ v=0, z=0}$ denote the projection parallel to $\partial_v$ onto the $(x,y)$-space. 
When $s\in \mathcal Op(0)$ goes to 0, the family $D'_1(s)$ accumulates 
to the $y$-axis in the $C^1 $-topology. 

Under the condition $\om_\phi(X_0)\neq 0$, that is $\cos_\phi(a^-)\neq 0$,
 Lemma \ref{lem:C1(s)} tell us that 
the family $C_1(s)\cap \{z=0\} $
accumulates  to the $x$-axis in the $C^1$-topology if and only if $s\,  \cos_\phi(a^-)$ goes to $0_+$. 
In particular, when $s\,  \cos_\phi(a^-)>0$ the projections $p_v(C_1(s))$ and $p_v(D'_1(s))$
intersect 
 in a unique point $b_1(s)$ and transversely. 
 If $s\,  \cos_\phi(a^-)$  is negative, then $C_1(s)\cap \{z=0\} $ is empty.\\

Denote by $c_1(s)\text{ and } d_1'(s)$ the only points in $C_1(s)\cap\sing{z=0}$ and in $D_1'(s)$ 
respectively such that $p_v(c_1(s))=b_1(s)=p_v(d_1'(s))$. Consider the real number
\begin{equation}\label{eq:v1(s)}
v_1(s):=v(c_1(s))-v(d_1'(s))\quad \text{ for every } s \text{ such that }s\,  \cos_\phi(a^-)\in \mathcal O_+.
\end{equation}
This function $v_1(s)$  depends smoothly on  $s$. Its derivative with respect to $s$
is  denoted by $\dot{v}_1(s)$.

\begin{remarque}\label{rem:g2}
By construction, $v_1(s)=0$ implies $c_1(s)=d_1'(s)$ which in turn implies 
 the existence of an orbit $\ell_s\in\Orb^{X_s}(p, p)$ passing through $c_1(s)$ such that $[\ell_s]=g^2$. This remark
 will be used when analysing the doubling phenomenon in Section \ref{section4}.
\end{remarque}

Lemma \ref{geometry-character} will show
 the kinematic 
 meaning of the character function $\chi$ at $X_0$.

\begin{lemme}\label{geometry-character} Let $(X_s)_s$ be a normalized 
positive crossing  path of $\S_g$ whose crossing velocity
\eqref{initial_velocity} is equal to $ +1$. If $\cos_\phi(a^-)\neq 0$ then
 the following relation holds: \begin{equation}\label{eq:CharGeom}
 \cos_\phi(a^-)\, \dot{v}_1(0)=\chi( X_0).
\end{equation}

\end{lemme}

\proof
Let us study the $v$-coordinate of $c_1(s)$ first. We notice that, if $\bar c_1(s)$ is another point 
of $C_1(s)$  depending smoothly on $s$ and such that $\bar c_1(0)=a^-=(\phi_0,0,-)$,
 we have the same velocity in $s= 0$:
 \begin{equation}\label{same-speed}
\frac{d}{ds} v\left(\bar c_1(s)\right)_{\vert s=0} =\frac{d}{ds} v\left(c_1(s)\right)_{\vert s=0}.
\end{equation}
Indeed, $C_1(s)$ accumulates $C^1$ to $\Si^-\cap\{z=0\}$ (Lemma \ref{lem:C1(s)}), then the difference $\dot c_1(0)-\dot{\bar c}_1(0)$ is a vector in $T_{a^-}\Si^-$.
We apply this remark to the point $\bar c_1(s):=C_1(s)\cap\{\phi=\phi_0\}$. 
Let $d_1(s)\in D_1(s)$  the lift of $\bar c_1(s)$ by $Desc^{-1}$.

Since $Desc$ preserves the $(\phi, r, \psi)$-coordinates, both paths 
$s\mapsto \bar c_1(s)$ and $s\mapsto d_1(s)$ have the same coordinates
 $(\phi(s)=\phi_0, r(s), \psi(s))$ when $s\neq 0$.
As $\bar c_1(s)\in\{\phi=\phi_0 \}$, the vector $\dot d_1(0)$ belongs to the $(n-i)$-plane 
which is the span of $\{\R\,\phi_0, T_{a^+}\Si^+\}$. 
Let ${\hat d}_1(0)$ 
be its projection to the line $\R\,{\phi_0}$ in the meridian disc $\{\psi=\psi_0\}$. 
Then,
\begin{equation}\label{eq:Info1}
\begin{array}{l}
 {\hat d}_1(0)=\rho\, \phi_0     
 \text{ where  } \rho={\ds\frac{d}{ds}r(s)_{\vert s=0} }.\\
\end{array}
 \end{equation} 
By definition of the $\phi$-latitude (Proposition \ref{prop:Latitudes})  
we have:
\begin{equation}
\langle \nu_\phi,{\hat d}_1(0)\rangle=\rho\langle \nu_\phi,\phi_0\rangle= \rho\, \cos_\phi (a^-).
 \end{equation}
  By definition, the hyperplane $\De^\phi$ is tangent in $a^+$ to   $D_1(0)$.
  Therefore, some calculus of Taylor expansion tells that  
 \begin{equation}
\langle \nu_\phi,{\hat d}_1(0)\rangle
=\frac{d}{ds} a^+(s)_{\vert s=0}=+1\,.
 \end{equation}
 We derive:
 \begin{equation}\label{eq:rho}
\rho=\frac{1}{ \cos_\phi(a^-)}\ .
\end{equation}

  Since $d_1(s) $ goes to $a^+=(-,0,\psi_0)$ as $s$ goes to $0$ and since the radial velocity is 
preserved by $Desc$, then we have:
\begin{equation}
\dot{\bar c}_1(0)= \rho\,\psi_0\in T_{a^-}\{\phi=\phi_0\}.
\end{equation}

Using again 
\eqref{eq:Trigo}, but relatively  
to  the preferred hyperplane $\Delta^{\psi}$  which defines the $\psi$-latitude
we obtain
\begin{equation}
\langle \nu_{\psi}, \dot{\bar c}_1(0)\rangle=\rho\,\cos_{\psi}(a^+)\,\nu_\psi.
\end{equation}
This together with the decomposition of $T_{a^-}\{z=0\}$ of \eqref{eq:decompTgt} says that 
there are  two vectors  $w_x\in T_{a^-}\Si^-$ and $ w_y\in \Delta^\psi$ such that
\begin{equation}\label{eq:c1'(0)}
\dot{\bar c}_1(0)=w_x+ \rho\,\om_{\psi}(a^+)\nu_{\psi} + w_y.
\end{equation}
 By\eqref{eq:c1'(0)}, we have 
$v(\dot{\bar c}_1(0))\, \partial_v^0
=\rho\,\om_{\psi}(a^+)\langle \p_v^0, \nu_{\psi}\rangle$. On the other hand, 
\eqref{eq:HolonDec} tells us that:
\begin{equation}\label{eq:HolomMod}
v(\nu_{\psi})=\frac{1}{\bar{\eta}}=\eta\,.
\end{equation}
Putting together (\ref{same-speed}),  \eqref{eq:rho}, \eqref{eq:c1'(0)} and \eqref{eq:HolomMod} we obtain:
\begin{equation}\label{eq:c1der}
v(\dot{c}_1(0))= v(\dot{\bar c}_1(0))=\rho\,\om_{\psi}(a^+)\, v(\nu_\psi)=\frac{\cos_{\psi}(a^+)}{\cos_{\phi}(a^-)}\,\eta\,.
\end{equation}\\

We come now to estimate the term $v(\dot{d'}_1(0))$. We apply Lemma \ref{velocity}
for comparing velocities associated with the holonomy $H_s$ and its inverse. From Formula (\ref{initial_velocity})
we derive that \hfill\break
$\frac {d}{ds} (v\circ H_s)(a^-)_{\vert s=0}= +1$. Then, the  inverse holonomy satisfies
\begin{equation}\label{second}
\frac {d}{ds} (v\circ H_s^{-1})(a^+)_{\vert s=0}= -1
\end{equation}
from which it is easily derived that $v(\dot{d'}_1(0))=-1$. Therefore:
\begin{equation}\label{eq:GrandV2}
\dot v_1(0)=\eta\frac{\cos_{\psi}(a^+)}{\cos_{\phi}(a^-)}+1
\end{equation}
which is a reformulation of the desired formula.\bull

Lemma \ref{lem:Passages} right below is the last tool that we need for proving Theorem 
\ref{thm:selfslideAlg}. It extracts the geometric information contained in Equation (\ref{eq:CharGeom}).
The setting is the same as in the previous lemma. We are only looking at normalized paths
$(X_s)_s$ which cross $\S_g$ positively at a point $X_0\in S_g^+$. 


\begin{lemme}\label{lem:Passages}${}$

\nd {\rm 1)} Suppose the $\phi$-latitude $ \cos_\phi(a^-)$ is positive  (and hence $\dot{v}_1(0)>0$).
Then, for $s\in \mathcal O_+$ there are sequences of non-empty $(i-1)$-discs $(D_k(s))_{k>1}$ and 
$(C_k(s))_{k>1}$ inductively defined  from the previous $D_1(s)$ and $C_1(s)$ by
\begin{equation}\label{defC_k}
\left\{ 
\begin{array}{l}
D_k(s):= H_s\left(C_{k-1}(s)\right) \cap \{z=1\}\\ 
C_k(s) := Desc\left( D_k(s)\right) \subset \partial^- \M_p
\end{array}
\right.
\end{equation}
 Moreover, as $s$ goes to 0,  the disc $C_k(s)$ tends to the North hemisphere  
 $\mathcal H^+(\Sigma^-)$
  in the $C^1$-topology, uniformly over every compact set of $\mathcal H^+(\Sigma^-)$.
  When $s\in \mathcal O_-$, both previous sequences are empty when $k>1$.
 \smallskip  
 
 \nd {\rm 2)} If $\dot{v}_1(0)$ and $ \cos_\phi(a^-)$ are negative, 
 then for $s\in \mathcal O_-$ the disc $C_2(s)$ is well defined as in \eqref{defC_k} 
 and the subsequent discs, $D_3(s), ...$, are empty. Moreover,
 $C_2(s)$ tends to $\mathcal H^+(\Si^-)$ in the   $C^1$-topology with the reversed orientation. 
 When  $s\in \mathcal O_+$, all discs in \eqref{defC_k}  are empty when $k>1$.
\end{lemme}


\proof 
1) When $s\in \mathcal O_-$, the disc $C_1(s)$ does not meet the tube $T$ around the 
homoclinic orbit $\ell$. Then $D_2(s)$ is empty and hence, all  further discs are so. 

Assume now that $s\in \mathcal O_+$. In that case, $C_1(s) $ goes to $\mathcal H^+(\Si^-)$
(Lemma \ref{lem:C1(s)})
and therefore meets the set $\{z=0\}$. Then,
the discs $D_2(s)$ and $C_2(s) $ defined in  \eqref{defC_k} are non-empty.
We are going to compute the position of $C_2(s)$ with respect to $D'_1(s)$ measured 
by some $v_2(s)$ in the direction
of the $v$-coordinate. 
We shall check the  positivity of  $\dot v_2(0)$ which will allow us to pursue the induction.

Recall the projection $\pi^{\psi_0}:\partial^+\M_p\to \{\psi=\psi_0\}$ and define the spherical annulus
$A:= (\pi^{\psi_0})^{-1}(\R\partial_v^1)$.
Consider the point $\tilde c_1(s)$ which is the transverse intersection
 $C_1(s)\cap H_s^{-1}(A)$. By projecting to 
the $v$-axis we find a function $v(\tilde c_1(s))$ which satisfies
\begin{equation}
\frac{d}{ds}v(\tilde c_1(s))_{\vert s=0}= \frac{d}{ds}v( c_1(s))_{\vert s=0}
\end{equation}
(Indeed, if $s$ goes to $0_+$, $\lim \tilde c_1(s)= \lim c_1(s)= a^-$).
Recall the definition of $d'_1(s)\in D_1'(s)$ from (\ref{eq:v1(s)}). 
Compute the derivative $V_2$ at $s=0$ of $v\left[H_s(\tilde c_1(s))\right]-v\left[H_s(d'_1(s))\right]$,
which is nothing but the velocity of the projection of $H_s(\tilde c_1(s))\in D_2(s)$ onto the $v$-axis of 
$\{z=1\}$ at $s=0$.
Using $\tilde c_1(0)= d'_1(0)=a^-$ and $dH_0(a^-)=Id$ in the coordinates of the tube $T$, we find:
\begin{equation}\label{V_2}
V_2= \dot v_1(0)
\end{equation}
which is positive by assumption. This $V_2$ will play the same r\^ole as the crossing velocity.

Since $V_2>0$, Lemma  \ref{lem:C1(s)} tells us that $C_2(s)$ meets $\{z=0\}$ when 
$s\in \mathcal O_+$. Therefore, we choose points $c_2(s)\in C_2(s)=Desc(D_2(s)) $ and
 $d'_2(s)\in D'_1(s)$ which forms the unique pair of points of the respective subsets 
 which have the same  $p_v$-projection. We define 
 \begin{equation}
 v_2(s):= v(c_2(s))-v(d'_2(s))
 \end{equation}
 The computation of $\dot v_2(0)$ is exactly the same except we have to replace 
 $V_1=1$ with $V_2$. The result is:
 \begin{equation}
 \dot v_2(0)= \eta\frac{ \cos_\psi(a^+)}{ \cos_\phi(a^-)} V_2+1.
 \end{equation}
 Here, 
 some discussion is needed according to the sign of $ \cos_\psi(a^+)$:
 \begin{itemize}
 \item[(i)] if $ \cos_\psi(a^+)$ is positive, then $\dot v_2(0)$ is larger than $ V_1=+1$. In that case
  the induction goes on with $V_k>V_{k-1}>\ldots> 1$.
  
\item[(ii)] if $ \cos_\psi(a^+)$ is negative, then $0<V_2=\dot v_1(0)<1$, where the last inequality comes from \eqref{eq:GrandV2}. Therefore, 
$\dot v_2(0)-1$  is the product of two numbers\footnote{One of them being $V_2$.} of opposite signs 
and whose absolute values are smaller than 1. Thus, $\dot v_2(0) $ belongs to $(0,1)$. Such a fact 
is preserved at each step of the induction.
 \end{itemize}
 The induction can be carried on.
 \smallskip
 
 \nd 2) Take $s\in \mathcal O_-$. The calculation yielding the equality (\ref{V_2}) still holds and tells us that
 $V_2$ is negative. Remark that $H_s(d'_1(s))\in \Si^+$. As $s<0$, one derives:
 $$
  v\left[H_s(\tilde c_1(s))\right]=v\left[H_s(\tilde c_1(s))\right]-v\left[H_s(d'_1(s))\right]>0.
 $$  Thus, Lemma \ref{lem:C1(s)} says that $C_2(s)$ tends to $\mathcal H^+(\Sigma^+)$
 in the $C^1$-topology. 
As $ \cos_\phi(a^-)<0$,  
$C_2(s)$ does not meet $\{z=0\}$ and the next discs are empty. 
 Concerning the orientation, we check that $D_2(s)$ tends to $-D_1(0)$ in $\partial ^+\M_p$. Then,
 $C_2(s) $ tends to $-\mathcal H^+(\Si^-)$. Finally, 
 the statement when $s\in \mathcal O_+$ is clear.
\bull

\begin{remarque}\label{rem-nonvanishing}
{\rm In the previous analysis, from Notation \ref{not:C1(s)} to Lemma \ref{lem:Passages}, we have 
given the lead role to the bottom of the Morse model, the attaching sphere $\Si^-$,
the perturbed holonomy and the map $Desc$. Here,
the non-vanishing of the $\phi$-latitude is required. 

One can make a similar analysis with the top of the Morse model, the co-sphere $\Si^+$, the inverse of 
perturbed holonomy and $Desc^{-1}$. There, the non-vanishing of the $\psi$-latitude is needed. But the 
the statements of the lemmas are analogous. As a consequence, if the proof of Theorem \ref{thm:selfslideAlg}
can be completed under the assumption $\om^\phi(X_0)\neq 0$, then it can also be completed 
when $\om^\psi(X_0)\neq 0$.
}
\end{remarque}

\subsection{\sc Proof of Theorem  \ref{thm:selfslideAlg} continued.\label{continued}}
We continue the proof which begins right after
 the statement of this theorem. After a series of 
equivalences, we are left to prove the case (1) of a positive crossing of the stratum $\S_g$
at a point $X_0\in \S_{g,L}$ where the character functon is positive. 
 We recall that 
 the statement of Theorem  \ref{thm:selfslideAlg} is preserved under the action of the groups 
 $\cG^\pm$ (see Notation \ref{cG}). Therefore, we may assume that $X_0\in \S_{g,L}$ is normalized.
 Moreover, as $\chi(X_0)\neq 0$, one of the extended $\phi$-latitude and  $\psi$-latitude is non-zero.
 By Remark \ref{rem-nonvanishing}, it is sufficient to complete the proof when $\om^\phi(X_0)\neq 0$.

The element $g\in \Pi_\al$ is thought of  as an arrow from the set of zeroes 
$Z(\al)$ into itself. Then $g$ determines its origin $p$ which is also its end point.
Recall that the Morse index of $p$ is $i$.
We look at any zero $q\in Z(\al)$ of Morse index $i-1$. We have to compute the change of 
$n(p,q)^X_L$ when $X$ changes from $X_{0_-}$ to $X_{0_+}$ in the given crossing path 
$\left(X_s\right)_s$.
It is useful to make some partition, adapted to $g$, of the set 
of connecting orbits from $p$ to $q$ for the gradient $X_{0_-}$.\\
 
 \nd {\sc Partition of the connecting orbits.} We may assume that each connecting orbit of 
 $X_{0_-}$ from $p$ to $q$ is the unique one in its homotopy class. In general, one would take the multiplicity into 
 account. Recall $\Ga_p^q$, the set of homotopy classes of paths from $p$ to $q$. The equivalence relation 
 defining the partition  of $\Ga_p^q$ is the following: $[\ga_0]\sim [\ga_1]$ 
 if and only if  the homotopy class of $\ga_1$  reads
  $[\ga_1]= g^k\cdot[\ga_0] \text{ with } k\in \Z$.

  Consider $[\ga]_{\sim}$, the $\sim$-class of a fixed connecting orbit $\ga$.
 Since the $\al$-lengths of connecting orbits are positive, we have  $u_\al([\ga'])<0$ for every $\ga'\in [\ga]_{\sim}$.
Therefore, as $u(g)<0$ there are only finitely many connecting orbits $\ga'\in [\ga]_{\sim}$ verifying
$u_\al([\ga'])\geq u_\al([\ga])$. Let $\ga_0$ 
be a connecting orbit in $[\ga]_{\sim}$ 
such that $u_\al([\ga_0])$
 is maximal. Then, any 
 element of $[\ga]_{\sim}$
 reads $g^k\cdot[\ga_0]$ for some 
 $k\geq 0$.\\
  
  \nd{\sc End of the proof.} By $\La_u$-linearity of the Morse-Novikov differential,
  without loss of generality we may assume that the above 
  partition has only one  $\sim$-class and that the maximal element $\ga$ is a positive connecting orbit 
  (with respect to the chosen orientations). 
   Let $b:= \ga\cap \Si^-$ and $\Delta_s$, $s=0_-$, be the connected component of $W^s(q,X_s)\cap\p^-\M_p$ containing $b$.
   After shrinking the parameter $\de$ of  $\M_p$ if necessary (see Subsection \ref{ssec:morse}), $\De_s$
  is an $(n-i)$-disc which  intersects $\Si^-$ transversely and only at $b$. 
  
   We are looking for the change formula 
   up to $\al$-lenght $L$ (for every $L>-u(g)$). Let $\left(X_s\right)_{s\in \cO p(0)}$ be a crossing path with $X(0)$ in
    $\S_g^+$. We do it first in the case where the $\phi$-latitude $\om(X(0))\neq 0$. 
   There are still four cases to consider
  where $a^-$ stands for $a^-(X_0)$ and $\mathcal H^\pm$ stand for $\mathcal H^\pm(\Si^-)$:
  \begin{enumerate}
   \item[(a.1)] The $\phi$-latitude $ \cos_\phi(a^-)$ is positive and $b$ belongs to  $\mathcal H^+$.
  \item[(a.2)] The $\phi$-latitude $ \cos_\phi(a^-)$ is positive and $b$ belongs to  $\mathcal H^-$.
   \item[(b.1)] The $\phi$-latitude $ \cos_\phi(a^-)$ is negative and $b$ belongs to 
   $\mathcal H^+$.
  \item[(b.2)] The $\phi$-latitude $ \cos_\phi(a^-)$ is negative and $b$ belongs to  $\mathcal H^-$.
 \end{enumerate}
  
  The proof consists of applying Lemma \ref{lem:Passages}. It is convenient to use the following 
  definiton. \\

\nd{\bf Definitions.}${}$\smallskip

\nd 1) {\it The positive (resp. negative) part of $W^u(p, X_0)$ is the union of all $X_0$-orbits passing 
through the positive (resp. negative) hemisphere $\mathcal H^+(\Si^-)$ (resp. $\mathcal H^-(\Si^-)$). It will be denoted 
by $W^u(p, X_0)^\pm$.}\smallskip

\nd 2) {\it For  a given $k>0$, we say that the unstable manifolds $W^u(p, X_s)$
  accumulate to $g^k\cdot W^u(p, X_0)^\pm$ when $s$ goes to $0_-$ (resp. $0_+$) if it is true when
  lifting to the universal cover, that is: 
  if $\tilde p$ (resp. $\widetilde X_s$)
   is a lift of $p$ (resp.  $X_s$), the unstable manifolds $W^u(\tilde p, \widetilde X_s)$
   accumulate to $W^u(g^k\tilde p, \widetilde X_0)^\pm$.}\\
 
 Here, it is worth noting that when a point lies in the accumulation set its whole
 {$X_0$-orbit}
  is also accumulated. As a consequence, Lemma \ref{lem:C1(s)} tells us that 
 $W^u(p, X_s)$ accumulate to \break
 $g\cdot W^u(p, X_0)^\pm$ in the $C^1$-topology
 when $s$ goes to $0_\pm$. Thanks to this topology,
 it makes sense to compare the orientations. The result is the following:  when $s\to 0_\pm$, then
 $W^u(p, X_s)$ accumulate to $\pm g\cdot W^u(p, X_0)^\pm$.  Accumulation to 
 $g^k\cdot W^u(p, X_0)^\pm$
 for some $k>1$ is dictated by Lemma \ref{lem:Passages} depending on the sign of 
 the $\phi$-latitude (knowing $\chi(X_0)>0$).
 We are now ready for proving (1) fromTheorem
 \ref{thm:selfslideAlg} in each above-enumerated case. 
 
 Let $\la_-(\ga)$ (resp. $\la_+(\ga)$) denote the element of the Novikov ring $\La_u$ which is 
  the contribution of  $[\ga]_{\sim}$ in $n(p,q)^{X_s}$ when $s<0$
  (resp. $s>0$). We have to check the next formula up to the given $L>0$ in each case (a.1) ... (b.2).
  \begin{equation}\label{change}
  \la_+(\ga)= (1+g+g^2+\cdots)\cdot \la_-(\ga)
  \end{equation}
  
  \nd {\bf Case (a.1).} When $s\to 0_-$, the oriented unstable manifolds $W^u(p, X_s)$  accumulate \break to $-g\cdot W^u(p,X_s)^-$ and nothing else. Therefore, 
 as $b\in \mathcal H^+$, we have 
$\la_-(\ga)= [\ga]$.
  
  When $s\to 0_+$, then $W^u(p, X_s)$ accumulate to $+g^k\cdot W^u(p,X_0)^+$ for every $k>0$ and will intersect $g^k\cdot\Delta_0$ transversely at a single point. Thus,  
  we have 
   $\la_+(\ga)=(1+g+g^2+\ldots)\cdot [\ga]$. 
   The change of $\la_\pm(\ga)$ from $s<0$ to $s>0$ 
  is really given by 
  Formula \eqref{change}.\\
   
   \nd {\bf Case (a.2).} As $b\in \mathcal H^-$ and taking into account the accumulation described 
   right above, we have:
   $ \la_-(\ga)= (1- g)\cdot [\ga]$ and $ \la_+(\ga)= [\ga]$. 
  Formula \eqref{change} is still fulfilled.\\
 
 \nd {\bf Case (b.1).} Here, the accumulation of $W^u(p, X_s)$ is
  dictated by part 2) of Lemma  \ref{lem:Passages} and the reason for Formula \eqref{change} 
   is more surprising than in the previous cases.
  When $s\to 0_-$,  the manifolds $W^u(p, X_s)$ accumulate to $-g\cdot W^u(p, X_0)^-$
  and to $-g^2\cdot W^u(p, X_0)^+$ and then nothing else. 
 When $s\to 0_+$, the manifolds $W^u(p, X_s)$ accumulate to
   $+g\cdot W^u(p, X_0)^+$ and nothing else.
   
   As $b\in \mathcal H^+$, we have $\la_-(\ga)= (1-g^2)\cdot [\ga]$ 
   and 
   $\la_+(\ga)= (1+g)\cdot [\ga]$. 
   Formula \eqref{change} is right since the identity 
   $(1+g+g^2+\cdots)(1-g^2) = 1+g$ holds in the Novikov ring.\\
   
  \nd {\bf Case (b.2).}  Accumulation is as right above. One derives that
  $\la_-(\ga)= (1-g)\cdot [\ga]$ and $\la_+(\ga)= [\ga]$. The desired formula is still satisfied.\\
  
  The proof of Theorem  \ref{thm:selfslideAlg} is now complete since only one $L$ is involved.
   \bull
   
   \begin{remarque}\label{virtual} One could ask what happens when there is no critical points
   $q$ of index $i(p)-1$. The answer is the following. The dichotomy $\S_g^+, \S_g^-$ still exists
   by the sign of the character function. Since
   the \emph{bifurcation factors} $(1+g+g^2+\cdots)$ and $(1+g)$ do not depend on $q$, one can associate them
   with each part of $\S_g\smallsetminus \S_g^0$ even if there is no $q$. In order to validate this association,
   it is sufficient to imagine a {\it virtual} zero $q^{virt}$ whose stable manifold intersects $\p_-\M_p$ along
    one generic  meridian at the beginning of a positive crossing path of $\S_g$. 
    The same proof as before tells us how  changes the number
     of virtual connecting orbits from $p$ to $q^{virt}$.
    \end{remarque}
   
   \subsection{\sc Proof of Theorem \ref{thm:selfslideSimplif}.} Here, the statement claims something
   to hold  \emph{for every 
    \break$L>-u(g)$} instead of \emph{for a given $L$}. In that case, it is natural that some genericity condition
   should be required. The condition in question---that is, a subset of $\S_g$---is the intersection 
   of all conditions: $X_0\in \S_{g,L}$  for $L\to+\infty$, each of them being the condition 
   which makes Theorem  \ref{thm:selfslideAlg} hold. A priori, this intersection could be empty.
   But thanks to Proposition \ref{g-residual}, this intersection is a residual set in $\S_g$ and we are done. \bull

   \section{\sc The doubling phenomenon. Proof of Theorem \ref{thm:Doubling}}\label{section4}

\subsection{\sc Notations and statement.}
In this section, we state and prove the refined version of Theorem \ref{thm:Doubling} which is given 
right  after specifying some definition and notations. It is about the local structure 
of  $\S_g^0$---the co-oriented 
locus in  $\S_g$ where the character function $\chi$ vanishes---in the complement of
$\S_g^{0,0}$, the latter being the locus where both of the  extended $\phi$-latitude  and $\psi$-latitude vanish.

\begin{defn} \label{decompositionS^0}
${}$

\nd {\rm 1)} Let $\R_+$ {\rm (}resp. $\R_-${\rm )} 
 be the set of positive {\rm (}resp. negative{\rm )}  real numbers.  
The open set $\S_g^{0,\pm}\subset \S_g^0$ is defined by the sign of the extended $\phi$-latitude, that is:
$X\in \S_g^{0,\pm} \Leftrightarrow \om_\phi(X)\in \R_\pm$.\smallskip

\nd {\rm 2)}  Let $X_{0,0}\in \S_g^0\prive \S^{0,0}_g$. 
Let  $\bigl(\cD(s,t):=X_{s,t}\bigr)$ 
be a \emph{germ} at $X_{0,0}$ of a two-parameter family in $\F_\al$, the space of adapted $\al$-gradients. 
This germ is said to be \emph{adapted} to the pair $(\S_g,\S_g^0)$ 
if the following conditions are fulfilled:
\begin{enumerate}
\item The one-parameter family $\bigl(\cD(0,t)\bigr)_{t\in \cO p(0)}$ is contained in $\S_g$, transverse to 
$\S_g^0$ and $\p_t\cD(0,0)$  is non-zero and points toward{s} $\S_g^+$. 
\item The partial derivative $\partial_s \cD(0,0)$ is transverse to $\S_g$ and points towards 
 its positive side.
\end{enumerate}
\end{defn}

\nd In particular, such a $\cD$ is transverse to $\S_g^0$. 

\begin{thm}\label{thm:DoublingRefined}
Let $\cD$ be 
a germ  of $2$-discs transverse  to $\S_g^0\prive \S_g^{0,0}$ and adapted to 
the pair $(\S_g,\S_g^0)$. Then $\cD$ intersects $\S_{g^2}$ transversely along an arc of $\S_{g^2}^+$. 
The trace on $\cD$ of the strata $\bparent{ \F_\al\,,\,\S_g\cup \S^+_{g^2}\, ,\, \S_g^{0,\pm}}$ 
is $C^1$-diffeomorphic  to 
  $$\bparent{\R^2,\, \R\times \{0\}\cup\{0\}\times \R_{\pm}\,,\, \{(0,0) \}}\,.$$
  Moreover, the natural co-orientation of $\S_{g^2}$ restricts to 
  the natural co-orientation of $\S_g^{0}$ in $\S_g$ or to its opposite depending upon 
  $\S_{g^2}$ approaches $\S_g^{0,+}$ or $\S_g^{0,-}$ respectively
  {\rm (Figure 2)}.
  
  Finally, if $X_{0,0}$ also fulfills the generic property  $\S_{g,\infty}$ {\rm (Definition \ref{almost})}
  then the germ $\cD$ does not meet $\S_{g^k}$ for $k>2$.
  \end{thm}
  
  Actually, the proof of Proposition \ref{g-residual} yields that the last property is generic in $\S^0_g$. Indeed,
 the  new constraint $\chi(X)=0$ only involves a compact domain of $W^u(p,X)$. 
   
   We first prove Theorem \ref{thm:DoublingRefined} for particular germs $\cD$ which we call
{\it  elementary}. Such a germ consists of  a one-parameter family of positive {\it normalized} crossing paths 
of $\S_g$ in the sense of Definition \ref{normalization} with some additional requirements. The definition of elementary 
path looks a bit strange, but it is inspired by a {\it toy model} of crossing $\S_g$ when all moving objects are
 affine subspaces in the coordinates of $\M_p$.

\subsection{\sc Elementary crossing path.}
Let $(X_s)_s$ be a normalized positive crossing path of  $\S_g$. 
After the normalization (Proposition \ref{conjugation}) we are still allowed 
  to prescribe more special dynamics of $X_s$;
   the {\it perturbed holonomy} will be specified near the respective homoclinic orbit $\ell$ of $X_0\in \S_g$   
   
  Let $a^\pm= \ell\cap \p_\pm\M_p$;  let $(\phi_0, 0,-)$ and $(-,0,\psi_0)$ be the respective muti-spherical coordinates
  of $a^-$ and $a^+$.
  Consider the spherical annulus  $\A_{\psi_0}:=\Si^-\times (0,1)\times \{\psi_0\}\subseteq \p^-\M_p$. 
  Assume the $\psi$-latitude of $a^+$ different from zero---by its very definition it is always the case when $X_0\in \S_g^{0,\pm}$.
  Therefore, whatever the perturbed holonomy $H_s$ along $\ell$
the inverse image $D'_1(s) $ of $\Si^+$ by 
$H_s$ 
 is transverse  to $\A_{\psi_0}$ for every $s$ close to 0. Call $b(s)$ the intersection point $D'_1(s)\cap \A_{\psi_0}$
when the intersection is non-empty; this is the case either when  $s< 0$ or $s>0$ 
depending on whether the $\psi$-latitude  $\cos_\psi(a^+)$ is \emph{negative} or \emph{positive}. 
By normalization, $b(s)$ belongs to the ray $\{(\phi_0,r,\psi_0)\mid r\geq 0\}$.
Below, we use Notation \ref{not:C1(s)}.

\begin{defn}\label{elementary}
The germ $\left(X_s\right)_{s}$ is said to be \emph{elementary} if it is normalized {\rm (}Definition \ref{normalization}
{\rm )}
and the following conditions are fulfilled.
\begin{enumerate}
\item The disc $D_1(s):= H_s(\Si^-)\cap\{z=1\}$ moves in the meridian disc $\{\psi=\psi_0\}$  
while remaining  parallel to 
the preferred  hyperplane $\De^\phi$ {\rm (\ref{eq:HypPhi})}. 
\item Let $a^+(s)$ be  the  intersection of $D_1(s)$ with the pole axis. For every $s$, 
\begin{equation} \label{elem-1}
 \p_s a^+(s) =1. 
 \end{equation}
\item For every $s$ close to 0 the  velocity of $b(s)$ 
is 
\begin{equation} \label{elem-2}
\p_s b(s) 
=-\frac{1}{\eta\,\cos_\psi(a^+)}\ .
 \end{equation}
 Here, $\eta$ stands for the holonomic factor of $X_0$ {\rm (Definition \ref{def:HolFactor})}.
\end{enumerate}
\end{defn}
 
 This definition makes sense only when the $\psi $-latitude of 
$a^+(X_0)$ is not $0$, that is, when $X_0$ does not lie on the $\phi$-axis $\S_g^\phi$ \eqref{axis}.  This
 is always the case when $X_0$ belongs to $\S_g^0\smallsetminus \S_g^{0,0}$.
 
 \begin{lemme}\label{lemme_elementary} Let 
$X_{0}\in \S_g\prive \S_g^\phi$ be an $\al$-gradient in normal 
 form. Then there exists a germ of elementary path $\left(X_s\right)$
 passing through $X_0$ and depending smoothly on $s$ in the $C^1$-topology. 
 \end{lemme}
 
 \proof Recall the tube $T$ with coordinates $(x,y,v,z)$ around the restricted  homoclinic orbit  $\underline{\ell}$ of $X_0$.
 The holonomy $H_0$ is defined on a neighborhood of $\{z=0\}$ in $\partial^-\M_p$
  and valued in a neighborhood of $\{z=1\}$ in $\partial^+\M_p$.
Since we are looking for a one-parameter perturbation $(X_s)$ of $X_0$ whose properties are readable in $T$ it is sufficient to describe it is near $T$.

 For $|s|$ small enough, the perturbed holonomy always reads 
  $H_s=H_0\circ K_s$ where $K_s$
 is a diffeomorphism of $\partial^-\M_p$ supported in its interior with $K_0= Id$.
 In order to satisfy conditions (\ref{elem-1}) and (\ref{elem-2}) of Definition \ref{elementary}, we first choose $a^+(s)$
 and $D_1(s)$ before choosing $K_s$.
For $a^+(s)$ we take the point in $\{\psi=\psi_0\}$ moving in the oriented pole axis  with velocity $+1$
 such that $a^+(0)=a^+$. For $D_1(s)$ we take the paralell disc to $D_1(0)$ passing through $a^+(s)$.
 Take $K_s$, smooth with respect to  $s$, such that:
 \begin{equation}\label{K_s}
 \left\{
 \begin{array}{l}
 K_s(a^-)= H_0^{-1}(a^+(s))\\
 K_s(\Si^-)
 =H_0^{-1}(D_1(s)) \text{ in } \cO p\{z=0\}
 \end{array}
 \right.
 \end{equation}
 Thus, the first two items are fulfilled. Note that by normalization of $X_0$ the point $K_s(a^-)$
 runs on a prescribed curve in the meridian $\{\phi=\phi_0\}$, namely the curve $H_0^{-1}(\R\nu_\phi)$. Its velocity 
 at time $s=0$ is the vector $\p_v^0$. By (\ref{eq:HolonDec}) we have
 \begin{equation}
 \langle \p_v^0, \nu_\psi\rangle= \frac 1\eta\,.
 \end{equation}
 
 We are now dealing with the last item. 
  We impose $D'(s)= K_s^{-1}(D'(0))$ to move in the meridian $\{\phi=\phi_0\}$; this is possible as the point
  $\kappa_s:= K_s^{-1}(a^-)$ already moves in this meridian by normalization. There are two more contraints: 
  the first one is $\p_s\kappa_s\vert_{s=0}= -\p_v^0$ by Lemma \ref{velocity}; the second one is (\ref{elem-2}). 
  This two constraints are compatible since $\langle \p_s b(s), \nu_\psi\rangle= -\frac 1\eta$.
  
  For having a one-parameter  family $(K_s)$ of diffeomorphisms of 
   $\p^-\M_p$ converging $C^1$ to identity when $s$ goes to 0,
  one has to choose conveniently $\p_sK_s$ at time $s=0$.  But is is easy to achieve since the velocity distribution
  is given along transverse submanifolds in the extended space $\left(\p^-\M_p\right)\times \cO p(0)$.
   \bull
   
   \begin{remarque} Note the great difference between the normalization process of a crossing path and the building
   of an elementary crossing path. The first one is achieved by an ambient $C^1$-conjugation; so, 
   it does not change the dynamics. The second one is a genuine bifurcation.
   \end{remarque}
 
 Clearly, Lemma \ref{lemme_elementary}  holds with parameters, for instance when the data is a one-parameter family 
 in $\S_g\prive \S_g^\phi$.
 Then, the next corollary follows.
 
 \begin{cor}\label{cor_elementary} Let 
$X_{0,0}\in \S_g^0\prive\S_g^{0,0}$ and let $\ga(t)=\left(X_{0,t}\right)_t$ be a  germ 
 of path in $\S_g$ passing through $X_{0,0}$ and crossing $\S^0_g$ transversely, such that 
 $\frac{\p\gamma}{\p t}(0)$ points towards $\S_g^+$. Then, there exists a
 two-parameter family $\cD=\left(X_{s,t}\right)$
 of pseudo-gradients of $\alpha$ adapted to $(\S_g,\S^0_g)$ such that,
 for every $t$ close to $0$, the path $s\mapsto X_{s,t}$ is elementary.
 Moreover, there  are such $X_{s,t}$ 
 which are smooth with respect to the parameters 
 in the $C^1$-topology.\\
 \end{cor}

\begin{defn}\label{two-elementary}
Let $\cD$ be a $2$-disc in $\F_\al$ transverse  to $\S_g^0\prive \S_g^{0,0}$ and adapted to the pair $(\S_g,\S_g^0)$. We say that $\cD$ is \emph{elementary} if it is made of a one-parameter family of elementary crossing paths as in {\rm Corollary \ref{cor_elementary}}.
\end{defn}

\nd{\bf Proof of Theorem \ref{thm:DoublingRefined}.}

First, we prove the theorem in the 
particular case where the transverse disc $\cD$
is elementary. 
Even in this particular case the proof is slightly different depending on where the base point
$X_{0,0}$ lies either in $\S_g^{0,-}$ or in $\S_g^{0,+}$. In each case, the proof has three items:
\begin{enumerate}
\item[\bf 1.] What is the trace  of $\S_{g^2}$ on $\cD$? 
Is there a non-empty trace of $\S_{g^k}$ for $k\neq 1\text{ or } 2$?
\item[\bf 2.] Is $\cD$ transverse to $\S_{g^2}$? How is the positive co-orientation of $\S_{g^2}$?
\item[\bf 3.]  Which 
part $\S_{g^2}^+$ or $\S_{g^2}^-$
is intersected  by $\cD$?\\
\end{enumerate}

\nd{\bf Case $X_{0,0}\in \S_g^{0,-}$.}  In other words, $a^-(X_{0,0})$ has a negative $\phi$-latitude.
\smallskip

\nd {\bf 1.} The pseudo-gradient  $X_{0,t}$ has a homoclinic orbit  $\ell_t$
based in $p$  and the $\phi$-latitude of
 $a^-(X_{0,t})$ lies in $[-1,0)$ for every $t\in \cO p(0)$. Denote by $\phi_t$ the spherical coordinate of 
 $a^-(X_{0,t})$. We use the tube $T$ around $\ell_0$ and its extremities: $\{z=0\}\subset \p^-\M_p$
 and $\{z=1\}\subset \p^+\M_p$.
 
 For simplicity, we specify  even more the path $\left(X_{0,t}\right)_t$ 
 by adding some assumptions (the discussion is similar with the other cases of latitudes by using
  other specifications\,\footnote{ If 
$\om_\phi(X_{0,0})=-1$, one makes $\frac{\p\eta}{\p t}(X_{0,t})>0$. Since 
$\om_\psi(X_{0,0})$ must be positive, $\frac{\p \chi}{\p t}(X_{0,t})>0$.}):
\begin{enumerate}
\item[(i)] The $\phi$-equator of $X_{0,t}$ is fixed and
the $\phi$-latitude $\cos_\phi(a^-(X_{0,0}))$ is not equal to $-1$.

\item[(ii)] The point $a^+(X_{0,t})= \Si^+\cap\ell_t$ and the $\psi$-equator of $X_{0,t}$ 
 are fixed.
\item[(iii)] The holonomic factor $\eta(X_{0,t})$ remains 
 constant and is denoted by $\eta$. 
\end{enumerate} 
Note
 that (i) allows one to take 
$(X_{0,t})_t$ positively transverse to $\S_g^0$ while satisfying (ii) and (iii). More precisely, one makes the $\phi$-coordinate  $\phi_t$ of $a^-(X_{0,t})$ vary on $t$ by increasing the $\phi$-latitude.

Denote  the spherical coordinate of $a^+(X_{0,t})$ by $\psi_0$, independent of $t$.
In this setting,  as the paths $s\mapsto X_{s,t}$ are elementary the discs $D_1(s,t)\subset \{z=1\}$
depend only on $s$ and are denoted by $D_1(s)$. For every $s\neq 0$, their 
images by the  descent map are 
discs $C_1(s)$ contained in the {\it spherical annulus} 
$\A_{\psi_0}:= \{(\phi,r,\psi_0)\mid \phi\in \Si^-, r\in [0,1]\}$. 
When $s$ goes to $0_-$,  by  Lemma \ref{lem:C1(s)}
the discs $C_1(s)$ accumulate to the negative hemisphere 
$\cH^-(\Si^-)$. 

Since $s\mapsto X_{s,t}$ is elementary and $\cos_\psi(\psi_0)>0$, the disc $D'_1(s,t)$, preimage in $\{z=0\}$ of $\Si^+$
by the respective perturbed holonomy, intersects $\A_{\psi_0}$ in one point $b(s,t)$ when $s\leq 0$ and nowhere when $s>0$, according to Definition \ref{elementary}. When $t$ is fixed, $b(s,t)$ moves on the ray 
 $\{(\phi_t, r,\psi_0)\mid r\geq 0\}$ and  its velocity 
is given by the formula in Definition \ref{elementary}.
According to Remark~\ref{rem:g2}, we have:
\begin{equation}\label{equationSg2}
X_{s,t}\in \S_{g^2}\quad \text{if and only if}\quad b(s,t)\in C_1(s).
\end{equation}

Denote by $c_1(s,t)$ the intersection point of 
$C_1(s)$ with the meridian disc $\{\phi= \phi_t\} $. When $t$
is fixed,  $c_1(s,t)$ also moves on the ray 
$\{(\phi_t, r,\psi_0) \mid r\geq 0\}$ and its radial velocity 
is the same as the one of its lift through $Desc$ in $D_1(s)\subset \p M_p^+$. Therefore,
\begin{equation}
\p_s c_1(s,t)= \frac {1}{\cos_\phi(\phi_t)}
\end{equation} 
As $X_{0,0}$ belongs to $\S_g^0$, that is $\chi(X_{0,0})=0$,
the 
curves $b(s,0)$ and $c_1(s,0)$, defined for $s<0$,
have the same radial velocities. Since both tend to 
$a^-(X_{0,0})$ on the same ray when $s$ goes to $0_-$,
we have $b(s,0)=c_1(s,0)$ for every $s$. Then, 
(\ref{equationSg2}) 
tells us that $X_{s,0}\in\S_{g^2}$ for every $s$ close to $0$ negatively.

For $t\neq 0$ and $s<0$, the radial velocities of $c_1(s,t)$ and $b(s,t)$ are distinct while 
their limits 
when $s$ goes to $0$ 
coincide. Therefore,  
(\ref{equationSg2}) tells us that $X_{s,t}$ never lies in $\S_{g^2}$ for $s<0$.

When $s>0$, the discs
$C_1(s) $ accumulate to the positive hemisphere $\cH^+(\Si^-)$. There is no chance 
for $C_1(s)$ to intersect 
$D'_1(s,t)$ which is 
far from any point in 
$\cH^+(\Si^-)$.

What about $\S_{g^k}$? If $k\leq 0$, we have $u(g^k)\geq 0$ and there is no homoclinic orbit in the homotopy class $g^k$. When $k>2$, we have to discuss the successive 
passages of the unstable manifold $W^u(p,X_{s,t})$ in $\p^-\M_p$, more precisely in 
$\{z=0\}$. 

By Lemma \ref{lem:Passages}, if $t>0$, that is $\chi(X_{0,t})>0$, and $s<0$ only the discs 
 $C_2(s,t)$ of the second passage 
 are  non-empty, but they accumulate to the positive hemisphere $\cH^+(\Si^-)$. 
 Therefore, no further passage could give rise to a  homoclinic orbit. When $s>0$, even the second passage does not exist.
 
 If $t<0$, one is able to see that there are infinitely many passages  in  $\{z=0\}$. But, by velocity considerations 
 $C_k(s,t)$ never meet $D'_1(s,t)$. We do not give more details here because this is similar
 to the symmetric case $X_{0,0}\in \S_g^{0,+}$ and $t>0$ where the analysis of velocities 
 will be completely achieved. Thus, the first item of case $X_{0,0}\in \S_g^{0,-}$ is proved.\\

 \nd{\bf 2.} The reason for transversality to $\S_{g^2}$ relies again on some 
computations of velocity.
 Define for $s\leq 0$:
 $$
 \de(s,t):= v\left(c_1(s,t)\right)-v\left(b_1(s,t) \right) \quad \text{and} \quad 
V(t):=\frac{\p\de(s,t)}{\p s} \vert_{s=0}\,.
$$
Although points are not the same, this velocity $V(t)$ at  $s=0$ is easily checked to be 
 the same as the velocity computed in Lemma \ref{geometry-character}. Then, for every $t$ close to $0$
 we have:
 \begin{equation}\label{V}
 V(t)= \eta\,\frac{\cos_{\psi}(\psi_0)}{\cos_\phi(\phi_t)} +1\quad \text{which implies} \quad \frac{dV(t)}{dt}<0\,.
 \end{equation}
 By definition of the character function, we have $V(0)=0$ which implies   $V(t)<0$ for $t>0$.
 Define  
$V(s,t):= \p_s\de (s,t)$.  By construction of $(X_{s,t})$, we have $V(s,0)=0$ for every $s<0$ close to 0.
 By (\ref{V}), the second partial derivative
 $\p^2_{t s}\de(s,t)$ is negative for every $(s,t)$ close to $(0,0)$ with $s\leq 0$ (here we use the 
 smoothness with respect to the parameters\footnote{The vector fields in a normalized crossing path are 
 not smooth with respect to the space variable. Their holonomy is $C^1$ only. Nevertheless, as the
  $C^1$-maps (of degree zero)  $\p^-\M_p\to \p^+\M_p$ form a Banach manifold it makes sense to 
  consider a smooth family of such holonomies.}). 
 By integrating 
 in the variable $s$ from $s_0<0$ to $0$ and noticing that $\de(0,t)=0$, we get:
\begin{equation}\label{pt}
 \frac{\p \de}{\p t} (s_0,t)= -\frac{\p}{\p t}  \left(\int_{s_0}^{0}\frac{\p \de}{\p  s}(s,t)\, ds\right) 
 = -\left(\int_{s_0}^{0} \p_{ts}^2\de(s,t)\, ds\right) >0.
 \end{equation}
 For $t=0$, this is exactly the transversality of $\cD$ to $\S_{g^2}$ at 
$X_{s_0,0}$. \\ 
 
  We are now looking at orientation. Take 
  $s_0<0$ such that  $b(s_0,0)$
 lies in $\{z=0\}$. It belongs to a homoclinic orbit $\ell'$ in the homotopy class $g^2$. There is a tube
 $T'$ around $\ell'$ with coordinates $(x',y', v', z') $. The $y'$-axis is contained in
 $D'_1(s_0,0)$ and is given a  co-orientation which follows from the co-orientation
  of $D'_1(0,0)$  by continuity.
 The $x'$-axis is contained in $C_1(s_0)$. Its  projection to the $x$-axis is orientation reversing
 (Lemma~\ref{lem:C1(s)}). Therefore:
 \begin{equation} \label{w}
 v(\p_{v'}) <0. 
 \end{equation}
  By (\ref{pt}) we have: 
${\ds \frac{\p}{\p t}[v\left(c_1(s_0,t)\right)-v\left(b_1(s_0,t) \right)]_{\vert_{t=0}} =   \frac{\p \de}{\p t} (s_0,0)>0}$.
  By replacing $v$ with $v'$ in the last inequality, we get: 
  $$ \frac{\p}{\p t}[v'\left(c_1(s,t)\right)-v'\left(b_1(s,t) \right)]_{\vert_{t=0}}<0\,.$$
  This translates the fact that $\p _t$ points to the negative side of $\S_{g^2}$ for $s<0$
  while for $s=0$, $\p _t$ defines the positive co-orientation of $\S_g^0$ in $\S_g$.\\

\nd{\bf 3.}  Let $L>0$. 
Consider a small  circle $\ga\subset \cD$ centered at the origin 
of the coordinates $(s,t)$ and turning positively with respect to the orientation given by these coordinates.
If the radius of $\ga$ is small enough\footnote{ The more $L$ is large, the more this radius has to be small.} 
 and if $X_{0,0}$ is generic,
$\ga$ avoids all codimension-one strata in $\mathcal F_\al$ except: 
\begin{itemize}
\item $\S_g$ which is crossed once in $\S_g^-$ positively, and once in $\S_g^+$ negatively,
\item $\S_{g^2}$ which is crossed once positively according to the above discussion.
\end{itemize}
As noted in Remark \ref{virtual} each crossed signed stratum is endowed with a \emph{bifurcation factor}.
The product of these factors should be equal to 1 up to $L$ in the Novikov ring after traversing $\ga$ once.
 The bifurcation factor of the  a small sub-arc of $\gamma$ crossing 
$\S_{g^2}$ is still unknown; call it $m(g)$.
This (commutative) product
is 
$$m(g)\cdot(1+g+g^2+\cdots)^{-1}\cdot (1+g)=1\quad{\rm(mod.}\ L). $$
Then, 
$m(g)=(1-g^2)^{-1} \ ({\rm mod.}\ L)$, that can only happen if 
the crossing of $\S_{g^2}$
takes place in  $\S_{g^2}^+$. The proof of  Theorem \ref{thm:DoublingRefined} is complete
for an elementary 2-disc in the case $X_{0,0}\in \S_g^{0,-}$.\\

\nd{\bf Case $X_{0,0}\in \S_g^{0,+}$.} In other words, $a^-(X_{0,0})$ has a positive $\phi$-latitude.
\smallskip

\nd{\bf 1.} The discussion is led in the same manner as in the previous case  
with  same notation. 
We only 
mention the main differences. Here,  $a^-(X_{0,0})$ belongs to the positive hemisphere of 
$\Si^-$ while $\psi_0$ belongs to the negative hemisphere of $\Si^+$. The discs $C_1(s)$ intersect
$\{z=0\}$ only when $s>0$. Therefore, for $s<0$ there is no chance for meeting  $\S_{g^k}$, 
for any $k\neq 0$.

By Lemma \ref{lem:Passages} there are 
infinitely many passages 
$C_k(s), \ k\geq 1, s>0$ of $W^u(p, X_{s,t})$ in $\p^-\M_p$ meeting $\{z=0\}$. Recall that 
the $(i-1)$-discs $C_k(s)$ 
do not depend on $t$.  The 
fact that $X_{s,t}$ belongs to 
$\S_{g^2}$ if and only if  
$s>0$ and $ t = 0$ is proved exactly  as 
in the previous case.

Then, we are left 
to show that for every $k>2$, $\S_{g^k}$ does not 
intersect $\cD$. Here it is important 
to think of $\cD$  as a germ because for a given representative this result is not true; when $k$
increases, the domain of the representative has to be restricted. 
Let $C_k(s,t)$ denote the $(i-1)$-disc in $\p^-\M_p$ 
corresponding to the $k$-th passage of the unstable manifold $W^u(p, X_{s,t})$
(see Lemma \ref{lem:Passages}); let $D'_1(s,t)$
denote the $(n-i-1)$-disc 
 corresponding to the first passage of the stable manifold $W^s(p, X_{s,t})$. 
Observe that
$\cD$ intersects $\S_{g^{k+1}}$ if and only if, for $(s,t)$ close to $(0,0)$, 
$C_k(s,t)$   intersects  $D'_1(s,t)$. This translates 
 in the next equation:
 \begin{equation}
 c_k(s,t)= d'_k(s,t)
 \end{equation}
where $c_k(s,t)$ and $d'_k(s,t)$ are 
the only two points of $\{z=0\}$ lying respectively on $C_k(s,t)$  and $D'_1(s,t)$
which have the same $(x,y)$-coordinates. Then, the above equation becomes:
 \begin{equation}\label{equality}
v\bigl( c_k(s,t)\bigr)= v\bigl(d'_k(s,t)\bigr)\,. 
 \end{equation}
When $s$ goes to $0$, these two points go to the same point 
$a^-(X_{s,t})\in\Si^-$. By computations done in the proof of Lemma \ref{lem:Passages},
we know that:
\begin{equation}
\frac{\p}{\p s}\left[v\bigl( c_k(s,t)\bigr)- v\bigl(d'_k(s,t)\bigr)\right]_{\vert_{s=0}}\neq 0\,.
 \end{equation}
It follows that, for $s$ close to $0$ (closeness depending on $t$), 
the equation (\ref{equality}) cannot be fulfilled. 

 The answer to questions {\bf 2} and {\bf 3} are exactly as in the case $X_{0,0}\in \S_g^{0,-}$. Then, Theorem \ref{thm:DoublingRefined}
is  proved for elementary 
2-discs as in Definition \ref{two-elementary}.\\

For finishing the proof of Theorem  \ref{thm:DoublingRefined} it is suitable to use  some $C^1$-topology. 
More precisely, we choose a system of finitely many
closed flow boxes $(B_j)_{j\in J}$ of $X_{0,0}$ whose end faces  $\p_\pm B_j$ are tangent to $\ker\al$
 and union
covers $M$ except a small open 
neighborhood $N$ of the zero set of $\al$.
 It is assumed that when slightly shrinking every $B_j$ to $B'_j$ tangentially to $\ker\al$, the union $\cup_j B'_j$
still covers $M\smallsetminus N$. Fix a  closed $C^0$-neighborhood $U$
of $X_{0,0}$ among the \emph{uniquely integrable} vector fields whose \emph{transverse
holonomy  is well defined}  for every $j\in J$ from $\p_-B'_j$ to $ \p_+B_j$ and of class $C^1$. 
This $U$ endowed with the $C^0$-topology of vector fields and the $C^1$-topology of holonomies $B'_j\to B_j,\, j\in J,$
may be thought of as a closed ball in a Banach manifold.

By Proposition \ref{conjugation}, there exists a neighborhood $V$ of $X_{0,0}$ in $\S_g^0$ 
such that the following properties are fulfilled for every $Y\in V$:
\begin{itemize}
\item[--] there exists a $C^1$-diffeomorphism $\Upsilon_Y$ of $M$ which 
carries $Y$ to a vector field in normal form, that is, $\left(\Upsilon_Y\right)_*Y$ is normalized;
\item[--] $\Upsilon_Y$  preverses the strata $\S_g$, $\S_g^0$ and $\S_g^{0,0}$ (Proposition \ref{invariance}).
\end{itemize}
It is easy to make $\Upsilon_Y$ depend continuously on $Y$ in the $C^1$-topology with the property that
$\Upsilon_Y= Id$ when $Y$ is already in normal form.

Let $\cD$ be an elementary 2-disc centered at $X_{0,0}$. If $Y\in V$ is close enough to $X_{0,0}$
and normalized one finds
an elementary 2-disc $\cD_Y$ centered at $Y$, depending $C^1$ on $Y$ and equal to $\cD$ if $Y= X_{0,0}$. 
If $Y\in V$ is not normalized
we still have a 2-disc centered at $Y$, namely
\begin{equation}
\cD_Y:=\left(\Upsilon_Y\right)_*^{-1}\cD_{\Upsilon_*Y}.
\end{equation}
This $\cD_Y$ is not elementary is the strict sense but it is conjugate to an elementary 2-disc in $U$. As 
$\Upsilon_Y$ preseves the stratification, in particular $\S_{g^2}$, the intersection of $\cD_Y$ with the 
different strata under consideration is the same as in the elementary case. Finally, we have a $C^1$-map
$$
F: V \times [-1,+1]^2\to U, \quad (Y,s,t)\mapsto F(Y,s,t)
$$
which meets $\S_{g^2}$ if and only if $t=0$ and $s\in \R_\pm $ depending upon $X_{0,0}\in \S_g^{0,\pm}$.
Moreover, the germ of $F$ at $(s,t)=(0,0)$ avoids all $\S_{g^k}$ for $k\neq 1,2$.

One checks that $span\{ \p_sF,\p_tF\}$ at $(s,t)=(0,0)$ is transverse to $\S_g^0$ in $U$. 
The Inverse Function Theorem is available and states that for  $V$ small enough 
 $F$ is a $C^1$-diffeomorphism 
onto its image $\mathcal N$, 
 an open set in $U$.
Therefore $\mathcal N$ has a product structure and a projection $P: \mathcal N\to [-1,+1]^2$ such that, for every $X\in \mathcal N$, the following equivalences hold:
\begin{equation}\label{equiv}
\left\{
\begin{array}{lcl}
 X \in \S_g & \Longleftrightarrow&  \bparent{s\circ P}(X)=0\\
 
 X \in \S_g ^0  &\Longleftrightarrow &   P(X)=(0,0)\\
  X \in \S_{g^2}& \Longleftrightarrow &   \bparent{t\circ P}(X)=0 \text{ and }
   \bparent{s\circ P}(X)\in \R_\pm 
  \text{ depending on }X_{0,0}\in \S_g^{0,\pm}.
  
\end{array}
\right.
\end{equation}
Let $\cD'$ be  any germ of two-parameter family centered in $X_{0,0}$ transverse to 
$\S_g^0\prive\S_g^{0,0}$ and contained in $\mathcal N$. 
Its projection $P\circ \cD'$ is submersive. The equivalences (\ref{equiv}) finish the proof of Theorem 
 \ref{thm:DoublingRefined}.
 
 \bull

\appendix \label{appendix}

\section{\\ Proof of the key fact \ref{deep}}
Let us recall the statement in question.

\begin{prop}\label{deep-appendix}
Let $X$ be an $\al$-gradient which is assumed Kupka-Smale. 
 Let $p$ and $q$ be two zeroes of $\al$
of respective Morse indices $k$ and $k-1$. Then, for every $L>0$ the number of connecting orbits from $p$
to $q$ whose $\al$-length is bounded by $L$ is finite.
\end{prop}

\proof The proof mainly consists of comparing the $\al$-length of any $X$-trajectory $\ga$
drawn on the unstable manifold $W^u(p, x)$ 
to the distance of its end points after lifting $\ga$ to the universal cover of $M$. By definition, the $\al$-length 
$\mathcal L(\ga)$ 
is additive with respect to any finite subivision of the considered trajectory $\ga$. \\

\nd {\sc First part.} 
For a trajectory $\ga$ descending from the top of a Morse model $\M(z)$ about any zero $z\in Z(\al)$ to the bottom 
of $\M(z)$ without 
getting out of it, $\mathcal L(\ga)$ is equal to the oscillation of any local primitive of $\al$ on $\M(z)$.
And hence, it does not depend on 
$\ga$. Without loss of generality, we may assume that this oscillation is the same for every zero of $\al$;
it is noted $h$. Therefore, given a trajectory $\ga$ of $\al$-length bounded by $L$,
the number $\kappa$ of segments traced on $\ga$
 by the compact union $\M:= \cup_{z\in Z(\al)}\M(z)$ fulfills
 \begin{equation}\label{number}
 \kappa\leq \frac Lh\, .
 \end{equation}

 In the complement of $\M$, that is away from the zeroes of $\al$, there is some positive constant $C$ such that 
 we have
 $\vert X(x)\vert \geq C$
for  every $x\in \M^*:= M\smallsetminus \M$.
 Denote by  $\la(\ga): = \int_\ga\vert\dot\ga\vert$ the length of a path $\ga$ and set 
 $\ga^*:= \ga\cap \M^*$.  Then, for every $X$-trajectory $\ga$ whose $\al$-length is bounded by $L$, we deduce 
 \begin{equation}\label{length*}
 C\la(\ell^*)= C\int_{\ell^*}\vert X(\ga(t))\, dt\vert  \leq  \int_{\ell^*}\vert X(\ga(t))\vert^2=\left\vert\int_{\ell^*}\al\right\vert\leq L\,. 
 \end{equation}
 Here the variable $t$ is the time of the flow of $X$.
 
Denote by $W_L(p)$ the union of the $X$-trajectories descending from $p$ whose $\al$-length is 
less than $L$. This is an open domain in the unstable manifold $W^u(p,X)$; it is homeomorphic to 
an open ball whose dimension is $k$. 
Let $\tilde M\mathop{\to}\limits^{\pi}M$
 be the universal cover of $M$ and let $\tilde p\in \pi^{-1}(p)$.  Let $\tilde X$ be the lift of $X$ to $\tilde M$. 
 It is a hyperbolic vector field and the unstable manifold $W^u(\tilde p, \tilde X)$ is the lift of $W^u(p,X)$
 through $\tilde p$. Moreover, its truncation $W^u_L(\tilde p, \tilde X)$ is the lift of $W^u_L(p, X)$.

 Let $\ell$ be an $X$-trajectory descending  from $p$ in $W^u_L(p, X)$, let  $e$ be its end point. Take its lift  $\tilde\ell$
 from $\tilde p$ and denote by $\tilde e$
  its end point. One looks at the subdivision $S$ of $\tilde\ell$ marked by its crossings with $\pi^{-1}(\M)$.
  From (\ref{number}), (\ref{length*}) and the triangular inequality applied to the vertices of $S$,
one deduces that the lifted distance satisfies
\begin{equation} \label{radius}
d(\tilde p, \tilde e)\leq \frac LC +\frac Lh\de=:R(L)
\end{equation} 
where $\de$ stands for the maximal diameter of $\M(z), \, z\in Z(\al)$. Therefore, we have an inclusion 
\begin{equation}
 W_L^u(\tilde p,\tilde X)\subset B(\tilde p, R(L))
\end{equation}
where $B(\tilde p, R)$ stands for the ball of $\tilde M $ about $\tilde p$ of radius $R$ and where $R(L)$
is the right hand side of (\ref{radius}). The consequence of these elementary estimations 
is that the closure $cl_L(\tilde p)$ of $\tilde W^u_L(\tilde p,\tilde X)$ is compact.
In particular, it contains finitely many lifts of the zero $q$ 
that we are interested in. Indeed these lifts cannot accumulate as their mutual distance is bounded from below.\\

\nd {\sc Second part.} The end of the proof uses the $KS$ assumption. Let $\tilde f$ be a global primitive
of $\pi^*\al$; it exists since $\tilde M$ is simply connected. The descending 
 gradient of $\tilde f$ with respect to the lifted metric
is equal to  $\tilde X$. 
The truncation of  $W^u(\tilde p, \tilde X)$
to the upper level set $\{\tilde f>\tilde f(\tilde p)-L\}$ is exactly the truncation $\tilde W^u_L(\tilde p,\tilde X)$.
The end point of a gradient line of $\al$-length $L$ belongs to the 
level set $\{\tilde f=\tilde f(\tilde p)-L\}$.

 Let us enumerate  $\tilde q_1, \ldots, \tilde q_m$ the lifts of $q$ which belongs to  $cl_L(\tilde p)$. 
 Now, we can argue as in usual Morse theory.
For $1\leq j\leq m$, consider the Morse model $\M(\tilde q_j)$ and the so-called co-sphere $\Si_j$, a sphere of 
dimension $(n-k)$ in the top of $\M(\tilde q_j)$. By the $KS$ assumption, the two following properties
hold for every $j=1, ..., m$:
\begin{enumerate}
\item the singular part (that is the frontier) of $cl_L(\tilde p)$ avoids $\Si_j$;
\item the regular part, that is $ W^u_L(\tilde p, \tilde X)$, is transverse to $\Si_j$.
\end{enumerate}
It is classical that the compactness of $cl_L(\tilde p)$ joined to these two properties implies the finitness of
$\Si_j\cap \tilde W_L(p)$. Therefore, there are finitely many orbits of $\tilde X$ descending from $\tilde p$
and ending at $\tilde q_j$ for every $j$. This is the desired finiteness. \bull


\vskip 1cm
\begin{thebibliography}{99}


\bibitem{farber} Farber M., {\it Topology of Closed One-Forms}, Math. Surveys and Monographs 108, Amer. Math. Soc., 2004.
\bibitem{h+w} Hatcher A. and Wagoner J., {\it Pseudo-isotopies of compact manifolds}, 
Ast\'erisque 6, Soc. Math. de France, Paris (1973). 
\bibitem{hutchings} Hutchings M., {\it Reidemeister torsion in generalized Morse theory}, 
Forum Math. 14 (2002), 209---244.
\bibitem{latour} Latour F., {\it Existence de {$1$}-formes ferm\'ees non singuli\`eres dans une classe de cohomologie de de {R}ham}, Inst. Hautes \'Etudes Sci. Publ. Math. 80 (1994), {135---194.}
{ \bibitem{flX}  Laudenbach F., {\it Transversalit\'e, Courants et Th\'eorie de Morse}, \'Editions de l'\'Ecole polytechnique, 
Ellipses, Paris, 2011.}
\bibitem{l-m} Laudenbach F., Moraga Ferr\'andiz C., {\it A geometric Morse-Novikov complex with infinite series 
coefficients}, C. R. Acad. Sci. Paris, Ser. I 356 (2018), 1222-1227. 
https://doi.org/10.1016/j.crma.2018.09.008
\bibitem{milnor} Milnor J.W., {\it Lectures on the $h$-cobordism theorem}, Princeton University Press, Princeton, NJ(1965).
\bibitem{moser} Moser J., {\it On the volume elements on a manifold,} Trans. Amer. Math. Soc. 
120 (1965), 286-294.
\bibitem{moragaTesis} Moraga Ferr\'andiz C., {\it Contribution \`a une th\'eorie de Morse-Novikov \`a param\`etre}, PhD thesis. Available at \href{http://hal.archives-ouvertes.fr/tel-00768575/}{http://hal.archives-ouvertes.fr/tel-00768575/},  Universit\'e de Nantes (2012).
\bibitem{moraga}  Moraga Ferr\'andiz C., {\it 
Elimination of extremal index zeroes from generic paths of closed 1-forms}, Math. Z.  278 (2014), 743---765.

\bibitem{novikov} Novikov S.P., {\it Multivalued functions and functionals. An analogue of the Morse theory}, Dokl. Akad. Nauk
SSSR 260(1), 31---35 (1981).
\bibitem{palis} Palis J., de Melo W., {\it  Geometric theory of dynamical systems : an introduction}, 
Springer-Verlag, New-York, 1982.

\bibitem{peixoto} Peixoto M.M., {\it On an approximation theorem of Kupka and Smale}, J. Differential Equations
3 (1966), 214---227.
{ \bibitem{shilnikov} Shilnikov L.P., {\it On the generation of a periodic motion from trajectories 
doubly asymptotic to an equilibrium state of saddle type}, Math. USSR Sbornik 6 (n$^\circ$3), 427---438. }
\bibitem{sikorav} Sikorav J.-C., {\it Points fixes de diff\'eomorphismes symplectiques, intersections de sous-vari\'et\'es lagrangiennes, et singularit\'es de 1-formes ferm\'ees}, th\`ese d'\'Etat, Paris 11, Orsay (1987). 

\end{thebibliography}
\end{document}